\documentclass{siamonline250211} %
\usepackage{lipsum}
\usepackage{amsfonts}
\usepackage{graphicx}
\usepackage{epstopdf}
\usepackage{algorithmic}
\ifpdf
  \DeclareGraphicsExtensions{.eps,.pdf,.png,.jpg}
\else
  \DeclareGraphicsExtensions{.eps}
\fi

\usepackage{array}
\usepackage{booktabs}
\usepackage[utf8]{inputenc}
\usepackage{amsmath,amssymb,bm}
\usepackage{mathrsfs}
\usepackage{enumitem}
\usepackage{overpic}
\usepackage{multirow}
\usepackage{nicefrac}
\usepackage{subfig}
\usepackage{float}
\usepackage{caption}
\usepackage{multicol}
\usepackage{dsfont}
\usepackage{bbm}
\usepackage{bm}
\usepackage{url}
\usepackage{verbatim}
\usepackage{algorithm}
\usepackage{xcolor}
\usepackage{cleveref}
\usepackage{hyperref}
\usepackage{pifont}
\hypersetup{
    colorlinks=true,
    linkcolor=blue,
    filecolor=blue,
    urlcolor=blue,
    citecolor=blue,
}

\usepackage{subfig} %[caption=false]
\usepackage{amsopn}

\graphicspath{
{images/}
}

\usepackage{enumitem}
\setlist[enumerate]{leftmargin=.5in}
\setlist[itemize]{leftmargin=.5in}

\newsiamremark{remark}{Remark}
\newsiamremark{hypothesis}{Hypothesis}
\crefname{hypothesis}{Hypothesis}{Hypotheses}
\newsiamthm{claim}{Claim}

\headers{Bijective Parameterization of \boldmath$n$-Manifolds}{T. Li, W.-W. Lin, Z.-H. Tan, X. Wan, S.-T. Yau, and J. Zhang}

\title{A Unified Framework for Bijective Joint Angle- and Volume-Preserving Parameterization of \boldmath$n$-Dimensional Manifolds
\thanks{Submitted to the editors DATE.
\funding{T. Li was partially supported by the National Natural Science Foundation of China (NSFC) No. 12371377 and the Jiangsu Provincial Scientific Research Center of Applied Mathematics under Grant No. BK20233002.}}}
\author{
Tiexiang Li\thanks{School of Mathematics and Shing-Tung Yau Center, Southeast University, Nanjing 211189, China; Shanghai Institute for Mathematics and Interdisciplinary Sciences (SIMIS), Shanghai 200433, China. (\email{txli@seu.edu.cn})}
\and Wen-Wei Lin\thanks{Shanghai Institute for Mathematics and Interdisciplinary Sciences (SIMIS), Shanghai 200433, China. (\email{wwlin@outlook.com})}
\and Zhong-Heng Tan\thanks{Department of Mathematics, The Chinese University of Hong Kong, Hong Kong, China. (\email{zhtan@math.cuhk.edu.hk})}
\and Xiao Wan\thanks{Center for Mathematics and Interdisciplinary Sciences, Fudan University, Shanghai 200433, China; Shanghai Institute for Mathematics and Interdisciplinary Sciences (SIMIS), Shanghai 200433, China. (\email{xwan25@m.fudan.edu.cn})}
\and Shing-Tung Yau\thanks{Yau Mathematical Sciences Center, Tsinghua University, Beijing 100084, China. (\email{styau@tsinghua.edu.cn})}
\and Junxin Zhang\thanks{Shanghai Institute for Mathematics and Interdisciplinary Sciences (SIMIS), Shanghai 200433, China; Research Institute of Intelligent Complex Systems, Fudan University, Shanghai 200433, China. (\email{jxzhang@simis.cn})}
}

\DeclareMathOperator{\diag}{diag}

\newcommand{\sgn}{\operatorname{sgn}}

\newcommand{\Det}{\operatorname{Det}}

\def\rmF{\mathrm{F}}
\def\rmR{\mathrm{R}}

\def\bfu{\mathbf{u}}
\def\bfv{\mathbf{v}}

\def\bfx{\mathbf{x}}
\def\bfy{\mathbf{y}}

\def\calM{\mathcal{M}}
\def\calN{\mathcal{N}}

\def\calX{\mathcal{X}}

\providecommand{\abs}[1]{\left\lvert#1\right\rvert}
\providecommand{\norm}[1]{\left\lVert#1\right\rVert}

\newcommand*{\dif}{\mathop{}\!\mathrm{d}}
\newcommand*{\re}{\mathbb{R}}

\newcommand*{\bff}{\mathbf{f}}

\newcommand*{\tr}{\operatorname{tr}}

\newcommand*{\SVol}{\operatorname{SVol}}
\newcommand*{\svol}{\operatorname{svol}}

\newcommand{\normmm}[1]{{\left\vert\kern-0.25ex\left\vert\kern-0.25ex\left\vert #1 
		\right\vert\kern-0.25ex\right\vert\kern-0.25ex\right\vert}}

\def\bbB{\mathbb{B}}

\def\bbS{\mathbb{S}}

\def\b1{\mathbbm{1}}

\def\ARAP{\mbox{\sc arap}}

\usepackage{arydshln}

\begin{document}

\maketitle

\begin{abstract}
    We propose a unified framework for bijective parameterizations of $n$-dimensional manifolds that jointly control angular and volumetric distortion through a weighted combination of conformal and volume-preserving energies. A logarithmic barrier based on signed simplex Jacobians prevents element inversion and degeneration during energy minimization. On oriented simplicial manifolds, the gradients of all three discrete energies admit a unified cotangent Laplacian-type representation that preserves the sparsity of the mesh connectivity. Based on these formulations, we develop a three-stage algorithm that constructs an initial parameterization, repairs foldings to restore strict feasibility, and subsequently minimizes the weighted distortion energy while maintaining local injectivity. The framework is formulated in a general setting, with particular attention to parameterizations onto the unit sphere $\mathbb S^n$ and the unit ball $\mathbb B^n$. {We derive sufficient conditions for discrete spherical maps to have and retain degree one, and apply classical topological criteria to establish global bijectivity.} Numerical experiments demonstrate that the proposed framework achieves simultaneous control of angular and volumetric distortion, while maintaining bijectivity and eliminating surface and volumetric folds. In a brain-tumor label-transfer study, the method improves tumor-label reconstruction compared with a non-bijective ablation.
\end{abstract}



\section{Introduction}
\label{sec:introduction}

Manifold parameterization constructs a bijective correspondence between a manifold and a canonical domain on which geometric data can be represented, compared, and processed. For an $n$-dimensional manifold $\calM$, the canonical target is typically an $n$-sphere, an $n$-ball, or an $n$-cube. Manifold parameterization is widely used as a fundamental tool in computer vision~\cite{MGBS21,LMKC13}, surface and volumetric registration~\cite{ABKC18,KCLM14,ZHTL25S,CPGP18,DZKC20}, remeshing~\cite{WCSS24,ECKC25,KSNL19}, morphing~\cite{MHWW17}, and medical image analysis and classification~\cite{WWCJ21,WWJW22,ZZHW25}. The quality of a parameterization depends on its global bijectivity and how well it preserves the source manifold's local geometry and volume.

Angle-preserving (conformal) parameterization controls angular and local shape distortion. A smooth map is conformal when its differential is, at every point, a composition of a rotation and a uniform scaling. Equivalently, all nonzero singular values of its tangent differential are equal. In two dimensions, conformal parameterization is supported by the rich theory of Riemann surfaces and has led to numerous numerical methods, including least-squares conformal maps~\cite{LBPS02}, discrete conformal parameterizations~\cite{MDMM02}, discrete Ricci flow~\cite{2008Discrete}, and boundary-first flattening~\cite{2017Boundary}. Conformal-energy minimization provides another efficient approach for constructing low-angle-distortion maps on simply and multiply connected surfaces~\cite{klyh:2021,LBPS02,MHWW19}.

The higher-dimensional problem is fundamentally more restrictive. By the
Liouville rigidity theorem, a conformal map between Euclidean domains in
dimension $n\geq3$ is the restriction of a M\"obius transformation~\cite[\S\,3.8.1]{GEFW13};
therefore, an exact conformal parameterization between prescribed
higher-dimensional manifolds generally does not exist. Higher-dimensional
methods consequently seek quasi-conformal maps by minimizing a distortion
functional rather than imposing exact conformality pointwise~\cite{IWT01,YTKC15,DPGP22}.

{Alongside angle preservation, volume preservation constitutes another fundamental consideration in manifold parameterization.} Its objective is to keep the local Jacobian determinant equal, or close, to a prescribed global volume ratio. For two-dimensional surfaces, density-equalizing maps and adaptive area-preserving parameterizations control area distortion through density diffusion or target-domain optimization~\cite{choi2022adaptive,choi2018density}. 
%谁谁提出了这个
Stretch-energy minimization was introduced as an energy-based approach to {area-preserving (equiareal)} surface parameterization~\cite{MHWW19} and was subsequently extended by the same line of work to volumetric manifolds~\cite{TMWH23,MHTL19} and, more generally, to $n$-dimensional volume-preserving parameterizations~\cite{ZHTL25}. Optimal mass transport provides another construction by matching prescribed volume distributions while minimizing the associated transportation cost~\cite{JWTL23,KSWC17,MHTM21}.

Angle- and volume-preserving objectives are, however, generally incompatible. Angle-preservation requires the singular values of the local differential to be equal, whereas volume preservation constrains only their product. An angle-preserving map may therefore exhibit substantial volume variation, while a volume-preserving map may contain severe angular and shape distortion. Moreover, volume preservation alone usually leaves a large family of admissible maps and does not select a unique or geometrically regular parameterization. It is thus natural to address the two objectives jointly: the angle-preserving component suppresses anisotropic deformation, whereas the volumetric component controls the local change of measure.

% 参考文献确认
Independently of distortion control, a valid parameterization must also prevent
element inversions and global overlaps. Existing strategies include weighted
Laplace--Beltrami eigen-projection for removing foldings in spherical conformal
maps~\cite{RLZW14}, orientation-preserving variational regularization for
three-dimensional image registration~\cite{DZKC20}, and scaffold-based reduction
of global injectivity to local flip-free constraints for planar
parameterization~\cite{QWWX23}. These methods address bijectivity in different
geometric settings and through different numerical mechanisms.

Recently, Liu and Yueh~\cite{LY25} proposed an energy-based distortion-balancing parameterization method for simply connected open triangular surfaces. Their formulation balances conformal and authalic distortions by imposing equality between the corresponding energy values. {Bijectivity is handled by switching to a mean-value Laplacian with positive weights under a prescribed convex boundary, a strategy specific to two-dimensional convex target domains.} For higher-dimensional Euclidean mappings, Lee, Lam, and Lui~\cite{YTKC15} developed a landmark-constrained variational model for bijective registration based on the conformality distortion $K_f(\bfx)=\frac{1}{n}\norm{Df(\bfx)}_{\rmF}^2/\det(Df(\bfx))^{2/n}$, with $K_f(\bfx)=+\infty$ when $\det(Df(\bfx))\leq0$. To address the computational demands of this nonlinear model, they introduced an auxiliary variable $R=Df$ and used ADMM to separate the global mapping update from elementwise positive-determinant subproblems.

Motivated by these developments, we propose a unified variational and numerical framework for bijective joint angle- and volume-preserving parameterizations of $n$-dimensional manifolds. At the continuous level, we introduce conformal and volume-preserving energies that measure deviations from local isotropy and uniform volume scaling, respectively. Both energies are nonnegative, and their zero-energy cases characterize angle- and volume-preserving mappings. Unlike the equality-constrained formulation in~\cite{LY25}, we combine them through a weighted objective, allowing a direct trade-off between angular and volumetric distortion. {We further prove that the joint energy controls the integrated as-rigid-as-possible (\ARAP) distortion, which measures the deviation of local deformations from rigid motions~\cite{LZXGG08}.} Furthermore, we introduce a logarithmic barrier energy based on signed simplex Jacobians to prevent element inversion and degeneration.

At the discrete level, the signed local volume ratio $J_\tau(f)$ records the volume change and orientation of each simplex $\tau$ and enters both the volume-preserving energy and the logarithmic barrier. Requiring $J_\tau(f)>0$ ensures that the affine map on each simplex is nondegenerate and orientation preserving, and hence injective. In addition to local injectivity, a degree-one topological condition is required to ensure global bijectivity~\cite{Hatcher2002,Lipman2014}.
We therefore derive sufficient approximation conditions for discrete spherical maps to have degree one and for subsequent refinement steps to preserve the degree-one property.

The numerical method follows the proposed variational formulation. A low-energy initial map is first constructed. If the initial map contains flipped simplices, an alternating direction method of multipliers (ADMM) scheme is employed to restore strict feasibility, where auxiliary positive Jacobian variables are introduced and the mapping subproblem is solved by a Riemannian preconditioned nonlinear conjugate-gradient method (PNCG). Once a strictly feasible map is obtained, a preconditioned Riemannian L-BFGS method is used to minimize the joint angle- and volume-preserving energy together with the logarithmic barrier.

These iterative optimization procedures require repeated evaluations of the corresponding energy gradients. To support efficient computation, based on~\cite{Crane19,ZHTL25}, we derive a unified cotangent Laplacian-type representation for the gradients of the conformal, volume-preserving, and barrier energies. This common representation preserves the sparsity induced by the mesh connectivity, enables simplex-local assembly, and allows the same computational framework and preconditioning strategy to be reused across the different energy terms.

The main contributions of this work are summarized as follows.
\begin{itemize}
\item \textbf{An \boldmath$n$-dimensional energy-based framework for angle and volume distortion control.} We develop a unified variational formulation for parameterizing $n$-dimensional manifolds with simultaneous control of angular and volumetric distortions. A weighted combination of conformal and volume-preserving energies enables flexible trade-offs between local anisotropy and volume change. We further prove that the joint energy controls integrated \ARAP\ distortion.

\item \textbf{A three-stage algorithm with global-bijectivity criteria.} We construct an initial parameterization, restore positive simplex Jacobians by ADMM--PNCG folding repair, and minimize distortion with a logarithmic barrier that prevents element inversion and degeneration. We derive sufficient conditions for discrete spherical maps to have degree one and to retain this degree during refinement, and apply classical global-bijectivity criteria to the spherical and ball settings.

\item \textbf{A unified cotangent Laplacian-type representation.} We derive discrete conformal, volume-preserving, and logarithmic barrier energies on oriented simplicial $n$-manifolds. Their gradients share the form $\nabla_{\mathbf f}E(f)=L(f)\mathbf f$, where $L(f)$ is a generalized cotangent Laplacian-type matrix. This structure supports simplex-local sparse assembly and allows the same solver and preconditioner framework to be reused across all energy terms.

\item \textbf{Medical-imaging application.} We combine bijective brain-to-ball and cube-to-ball maps to transfer tumor labels between patient-specific brain domains and a canonical cube. Experiments show improved tumor-label reconstruction compared with a non-bijective ablation.
\end{itemize}

The remainder of the paper is organized as follows. \Cref{sec:ContinuousEnergy} introduces the continuous conformal, volume-preserving, and joint angle- and volume-preserving energies, characterizes their zero-energy mappings, and establishes the \ARAP\ distortion bound for the joint energy. \Cref{sec:DiscreteEnergy} develops the corresponding discrete energies and establishes their unified cotangent Laplacian-type gradient formulas. \Cref{sec:alg} presents the three-stage algorithms for spherical and ball parameterizations, together with the analysis of degree preservation and global bijectivity. The final sections report numerical experiments and a medical-imaging application concerning brain-tumor label transfer.

\section{\boldmath$n$-dimensional metric-preserving energies}\label{sec:ContinuousEnergy}

In this section, we develop a continuous energy and variational framework for metric-preserving parameterizations. We first introduce the $n$-dimensional conformal and volume-preserving energies for smooth maps between manifolds and characterize their zero-energy mappings. These energy functionals provide quantitative measures of angular and volumetric distortions. Since exact conformality is generally too restrictive in higher dimensions, while volume preservation alone usually leaves considerable freedom and may lead to large angular distortions, we combine the two energies into a joint formulation that promotes locally isotropic deformation while controlling local volume change. We further incorporate local nondegeneracy requirements to exclude collapsed or folded parameterizations. The resulting constrained variational model lays the foundation for the discrete energies and numerical methods developed in the subsequent sections.

Let $\calM$ be a smooth compact $n$-dimensional manifold embedded in $\mathbb{R}^m$ ($2\leq n\leq m$). For each $\bfx\in\calM$, there exist an open set $\Omega_{\bfx}\subseteq\mathbb{R}^n$, a neighborhood $U_{\bfx}\subseteq\calM$ of $\bfx$, and a smooth local parametrization $h:\Omega_{\bfx}\rightarrow U_{\bfx}$. Its Jacobian $J_h(\bfu)$ has full column rank for every $\bfu\in\Omega_{\bfx}$.
Let $f:\calM\rightarrow\mathcal{N}$ be an orientation-preserving and bijective smooth map, where $\mathcal{N}$ is an $n$-dimensional target manifold embedded in a Euclidean space. On each local parametrization $h$, the tangent gradient of $f$ on $\calM$ is defined as $\nabla_{\calM}f:=J_{f\circ h}J_h^\dagger=J_{f\circ h}(J_h^\top J_h)^{-1}J_h^\top$. This definition is independent of the choice of local parametrization. Let $\sigma_1\geq\cdots\geq\sigma_n>0$ be the singular values of $\nabla_{\calM}f$, and let $\Det(\cdot)$ denote the pseudo-determinant, i.e., the product of all nonzero singular values. Then $\|\nabla_{\calM}f\|_{\rmF}^2=\sum_{i=1}^n\sigma_i^2$, $\Det(\nabla_{\calM}f)=\prod_{i=1}^n\sigma_i$.

\subsection{Conformal energy}

The \textbf{Dirichlet energy} of $f:\calM\to\calN$ is defined by
\begin{equation}\label{eq:Dirichlet_Energy}
    E_D(f)=\frac{1}{n}\int_{\calM}\norm{\nabla_{\calM}f}_{\rmF}^2\dif\mu=\frac{1}{n}\int_{\calM}\sum_{i=1}^n\sigma_i^2\dif\mu,
\end{equation}
where $\dif\mu$ is the intrinsic $n$-dimensional volume element on an $n$-manifold $\calM$.

For $n=2$, the classical conformal energy in \cite{klyh:2021} is defined as the difference between the Dirichlet energy and the area of the image. On the other hand, the pointwise quantity $K_f(\bfx):=\frac{1}{n}\norm{\nabla_{\calM}f}_{\rmF}^2/\left(\Det(\nabla_{\calM}f)\right)^{2/n}$ is the generalized conformality (angle) distortion used in the theory and computation of $n$-dimensional quasi-conformal maps; see, for example,~\cite{YTKC15}. Motivated by the two-dimensional conformal energy and the higher-dimensional conformality (angle) distortion, we define the $n$-dimensional \textbf{conformal energy}
\begin{equation}\label{eq:Conformal_Energy}
E_C(f):=E_D(f)-E_A(f),
\end{equation}
where
\begin{equation}\label{eq:Area_energy}
E_A(f):=\int_{\calM}\left(\Det(\nabla_{\calM}f)\right)^{2/n}\dif\mu=\int_{\calM}\left(\prod_{i=1}^n\sigma_i\right)^{2/n}\dif\mu.
\end{equation}
When $n=2$, $E_A(f)$ equals the area of $f(\calM)$, and hence \cref{eq:Conformal_Energy} reduces to the classical conformal energy in \cite{klyh:2021}. For this reason, we refer to $E_A(f)$ as the \textbf{area energy}. The following result is a direct consequence of the generalized conformality distortion in \cite[Definition~1]{YTKC15} and extends the corresponding property of the two-dimensional conformal energy.
\begin{theorem}[\cite{YTKC15}]\label{thm:conformal_energy}
The conformal energy satisfies $E_C(f)\geq 0$, and equality holds if and only if $f$ is conformal at every point $\bfx\in\calM$, i.e., $\sigma_1=\sigma_2=\cdots=\sigma_n$.
\end{theorem}

This theorem shows that $E_C(f)$ measures the deviation of $f$ from conformality, with zero energy attained only by exact conformal parameterizations. Since such parameterizations generally do not exist, minimizing $E_C(f)$ seeks a mapping with the smallest conformality (angle) distortion within the admissible class.

\subsection{Volume-preserving energy}
Let $\nu$ denote the standard intrinsic $n$-dimensional volume measure on the target manifold $\mathcal{N}$. The map $f$ is said to be volume-preserving if
\begin{equation}\label{eq:vol_pres}
    \frac{\nu(A)}{\nu(f(\calM))}=\frac{\mu(f^{-1}(A))}{\mu(\calM)}
\end{equation}
for every Borel-measurable set $A\subseteq f(\calM)$. By the change-of-variables formula on embedded manifolds, \cref{eq:vol_pres} is equivalent to $\Det(\nabla_{\calM}f)=\rmR_f$ at each point $\bfx\in\calM$, where $\rmR_f:=\nu(f(\calM))/\mu(\calM)$ is the global volume ratio induced by $f$, namely, the ratio between the intrinsic $n$-dimensional volume of the image $f(\calM)$ and that of the source manifold $\calM$. Accordingly, we define the \textbf{volume-preserving energy} by
\begin{equation}\label{eq:EVcontinuous}
    E_V(f)=\int_{\calM}\left(\Det(\nabla_{\calM}f)-\rmR_f\right)^2\dif\mu.
\end{equation}
Tan et al.~\cite{ZHTL25} defined the \textbf{volumetric stretch energy} of $f$ as
\begin{equation}\label{eq:EVScontinuous}
    E_{V_s}(f)=\int_{\calM}\left(\Det(\nabla_{\calM}f)\right)^2\dif\mu.
\end{equation}
By the change-of-variables formula, the energy $E_V(f)$ in \cref{eq:EVcontinuous} reduces to
\begin{equation}\label{eq:EVIDcontinuous}\begin{aligned}
    E_V(f)&=E_{V_s}(f)-2\rmR_f\int_{\calM}\Det(\nabla_{\calM}f)\dif\mu+\rmR_f^2\cdot\mu(\calM)\\
    &=E_{V_s}(f)-2\frac{[\nu(f(\calM))]^2}{\mu(\calM)}+\frac{[\nu(f(\calM))]^2}{\mu(\calM)}=E_{V_s}(f)-\frac{[\nu(f(\calM))]^2}{\mu(\calM)}.
\end{aligned}\end{equation}
We refer to $E_{V_r}(f):={[\nu(f(\calM))]^2}/{\mu(\calM)}$ as the \textbf{volumetric ratio energy}. If $f(\calM)=\calN$, then $E_{V_r}\equiv E_{V_r}(f)$ is a constant, and the volume-preserving energy is equivalent to $E_{V_s}(f)$ up to the additive constant $E_{V_r}$. The following result is obvious in view of the definition \cref{eq:EVcontinuous} of the volume-preserving energy.
\begin{theorem}\label{thm:volumetric_energy}
    The volume-preserving energy satisfies $E_V(f)\geq0$, with equality if and only if $f$ is volume-preserving up to the global volume ratio $\rmR_f$.
\end{theorem}

Therefore, the volume-preserving parameterizations are precisely the global minimizers of $E_V(f)$ with minimum value zero.

\subsection{Bijective joint angle- and volume-preserving parameterization}

For $n\geq3$, exact conformality imposes a highly restrictive local isotropy condition, and an angle-preserving parameterization between prescribed source $\calM$ and target manifold $\calN$ generally does not exist. In contrast, volume preservation alone typically admits a large family of feasible parameterizations and therefore leaves considerable freedom in the local deformation. This freedom may lead to severe angular and shape distortions, even when the local volumes are perfectly preserved. It is therefore natural to impose, to an appropriate extent, the more rigid angle-preservation requirement within the volume-preserving formulation. The angular term reduces the otherwise excessive degrees of freedom and favors locally isotropic deformations, while the volumetric term controls the local volume change. Their combination therefore balances local isotropy and volume preservation within a single variational formulation.

We define the \textbf{joint angle- and volume-preserving energy} by
\begin{equation}\label{eq:Joint_Energy}
    E_{\beta}(f)=\beta E_C(f)+(1-\beta)E_V(f),\qquad 0\leq\beta\leq1.
\end{equation}
The parameter $\beta$ controls the relative weighting of angle and volume preservation. The cases $\beta=1$ and $\beta=0$ correspond to pure angle-preserving and pure volume-preserving energy, respectively.

For $\beta\in(0,1)$, the ideal condition $E_\beta(f)=0$,
although generally unattainable for prescribed source and
target manifolds, implies the singular values of the tangent gradient satisfy
$\sigma_1=\cdots=\sigma_n=\rmR_f^{1\over n}$ at each $\bfx\in\calM$.
A quantitative version of this observation relates the joint energy
to as-rigid-as-possible (\ARAP) distortion.
In the surface parameterization method of Liu et al.~\cite{LZXGG08}, \ARAP\ distortion is measured by the squared Frobenius distance from
each triangle differential to its closest rotation. For an
orientation-preserving differential expressed in oriented orthonormal
tangent frames, this distance is $\sum_{i=1}^n(\sigma_i-1)^2$.
To account for the global volume ratio, define the \ARAP\ energy at this prescribed scale by
\begin{equation}\label{eq:scaled_ARAP_energy}
    E_{\ARAP}(f)
    :=\int_{\calM}\sum_{i=1}^n
    \left(\sigma_i-\rmR_f^{1\over n}\right)^2\dif\mu.
\end{equation}
For bijections between prescribed source and target manifolds,
$\rmR_f=\nu(\mathcal N)/\mu(\calM)$ is fixed. When $\rmR_f=1$,
\cref{eq:scaled_ARAP_energy} is the usual continuous \ARAP\ energy.
The following estimate bounds this distortion in terms of
$E_C(f)$ and $E_V(f)$.

\begin{theorem}\label{prop:ARAP_control}
The \ARAP\ energy satisfies
\begin{equation}\label{eq:ARAP_moment_bound}
    E_{\ARAP}(f)
    \leq nE_C(f)+n\mu(\calM)\rmR_f^{2\over n}
    \Phi_n\left(\frac{E_V(f)}{\mu(\calM)\rmR_f^2}\right),
\end{equation}
where, for $v\geq0$, $\Phi_n(v)=1+[\rho_n(v)]^{n-2}-2[\rho_n(v)]^{n-1}$, $\rho_n(v)=\max\left\{(1+v)^{-{1\over n}},\frac{n-2}{2(n-1)}\right\}$. The function $\Phi_n$ satisfies $0\leq\Phi_n(v)\leq v$.
Consequently, for every $\beta\in(0,1)$,
\begin{equation}\label{eq:ARAP_joint_bound}
    E_{\ARAP}(f)
    \leq n\max\left\{\frac1\beta,
    \frac{\rmR_f^{{2\over n}-2}}{1-\beta}\right\}E_\beta(f).
\end{equation}
\end{theorem}

\begin{proof}
Write $V=\mu(\calM)$ and $J=\Det(\nabla_{\calM}f)=\prod_{i=1}^n\sigma_i$ for simplicity. The arithmetic--geometric mean inequality gives
$\sum_{i=1}^n\sigma_i\geq nJ^{1\over n}$, and hence
\begin{equation}\label{eq:localAMGM}
    \sum_{i=1}^n(\sigma_i-\rmR_f^{1\over n})^2
    \leq \sum_{i=1}^n\sigma_i^2+n\rmR_f^{2\over n}-2n\rmR_f^{1\over n}J^{1\over n}=\sum_{i=1}^n\sigma_i^2-nJ^{2\over n}
       +n(J^{1\over n}-\rmR_f^{1\over n})^2.
\end{equation}
Integrating \cref{eq:localAMGM} over $\calM$ gives
\begin{equation}\label{eq:ARAP_shape_scale_split}
    E_{\ARAP}(f)
    \leq nE_C(f)+n\int_{\calM}(J^{1\over n}-\rmR_f^{1\over n})^2\dif\mu.
\end{equation}

To estimate the remaining integral, we set $v=E_V(f)/(V\rmR_f^2)$ and obtain that
\begin{equation}\label{eq:J-moment}
    \int_{\calM}J\dif\mu=\nu(f(\calM))=V\rmR_f,\quad\int_{\calM}J^2\dif\mu=E_{V_s}(f)=E_V(f)+V\mathrm{R}_f^2=(1+v)V\mathrm{R}_f^2.
\end{equation}
For $x>0$, define $\varphi_n(x)=x^{-1}(x^{2\over n}-2x^{1\over n})$, which satisfies $(x^{1\over n}-1)^2=1+x\varphi_n(x)$.
Direct differentiation gives
\[\varphi_n'(x)=\frac{x^{{1\over n}-2}}{n}\left[2(n-1)-(n-2)x^{1\over n}\right],\quad\varphi_n''(x)=\frac{2(n-1)}{n^2}x^{{1\over n}-3}\left[(n-2)x^{1\over n}-(2n-1)\right].\]
Thus $\varphi_n$ is increasing and concave on $(0,L_n)$, where $L_n=\infty$ for $n=2$ and $L_n=\left(2{n-1\over n-2}\right)^n$ for $n\geq 3$. Let $\widehat{\varphi}_n$ be the truncation of $\varphi_n$ at the level $L_n$, i.e., $\widehat{\varphi}_n(x)=\varphi_n(\min\{x,L_n\})$, which is nondecreasing and concave on $(0,\infty)$, with $\widehat \varphi_n\geq \varphi_n$. Applying Jensen's inequality with weight $J/(V\rmR_f)$ and using the equalities in \cref{eq:J-moment}, we obtain that
\begin{equation}\label{eq:PfYUNDAO}
    \begin{aligned}
        &\int_{\calM}(J^{1\over n}-\mathrm{R}_f^{1\over n})^2\dif\mu=\mathrm{R}_f^{2\over n}\int_{\calM}\left(1+{J\over \mathrm{R}_f}\varphi_n\left({J\over\mathrm{R}_f}\right)\right)\dif\mu\\
        =\,&V\mathrm{R}_f^{2\over n}\left(1+\int_{\calM}{J\over V\mathrm{R}_f}\varphi_n\left({J\over\mathrm{R}_f}\right)\dif\mu\right)\leq V\mathrm{R}_f^{2\over n}\left(1+\int_{\calM}{J\over V\mathrm{R}_f}\widehat{\varphi}_n\left({J\over\mathrm{R}_f}\right)\dif\mu\right)\\
        \leq\,&V\mathrm{R}_f^{2\over n}\left(1+\widehat{\varphi}_n\left(\int_{\calM}{J^2\over V\mathrm{R}_f
        ^2}\dif\mu\right)\right)=V\mathrm{R}_f^{2\over n}(1+\widehat{\varphi}_n(1+v))=V\mathrm{R}_f^{2\over n}\Phi_n(v).
    \end{aligned}
\end{equation}
The last identity follows from $\rho_n(v)=\min\{1+v,L_n\}^{-{1\over n}}$. Then \cref{eq:ARAP_moment_bound} is proved by substituting \cref{eq:PfYUNDAO} into \cref{eq:ARAP_shape_scale_split}. Finally, $\widehat \varphi_n(1)=-1$ and $\widehat \varphi_n'(1)=1$. Monotonicity and concavity therefore give $0\leq\Phi_n(v)=1+\widehat \varphi_n(1+v)\leq v$. Hence
\begin{equation*}
    E_{\ARAP}(f)\leq nE_C(f)+n\rmR_f^{{2\over n}-2}E_V(f)\leq n\max\left\{\frac1\beta,\frac{\rmR_f^{{2\over n}-2}}{1-\beta}\right\}E_\beta(f),
\end{equation*}
which completes the proof.
\end{proof}

Thus, for each fixed $0<\beta<1$, small joint angle- and volume-preserving energy guarantees small integrated \ARAP\ distortion at the prescribed global scale. When $\mathrm{R}_f=1$, from \cref{eq:ARAP_joint_bound} we deduce that $E_{\ARAP}(f)\leq 2nE_{1\over 2}(f)$.

\begin{remark}\label{rem:ARAP_volume}
    The variance bound for the arithmetic--geometric mean difference~\cite{Aldaz2012} gives
    \begin{align*}
        &\frac1n\sum_{i=1}^n\sigma_i^2-\left(\prod_{i=1}^n\sigma_i\right)^{2/n}
        \leq\sum_{i=1}^n\left(\sigma_i-\frac1n\sum_{j=1}^n\sigma_j\right)^2\\
        =\,&\sum_{i=1}^n(\sigma_i-\mathrm{R}_f^{1/n})^2-n\left(\frac1n\sum_{j=1}^n\sigma_j-\mathrm{R}_f^{1/n}\right)^2\leq\sum_{i=1}^n(\sigma_i-\mathrm{R}_f^{1/n})^2.
    \end{align*}
    Integration over $\calM$ yields $E_C(f)\leq E_{\ARAP}(f)$. Thus, sufficiently small \ARAP\ energy guarantees small integrated conformal distortion. However, small integrated \ARAP\ energy can coexist with severe local volume changes. To see this, let $\mathbb{B}^n$ denote the closed unit ball in $\re^n$. For $M>1$ and $0<\varepsilon<1$, choose $g_{M,\varepsilon}\in C^\infty([0,1])$ with $0<g_{M,\varepsilon}'\leq M$ such that
    \begin{equation*}
        g_{M,\varepsilon}(r)=
        \begin{cases}
            Mr,&0\leq r\leq \dfrac{\varepsilon}{4M},\\[4pt]
            r,&\dfrac{\varepsilon}{2}\leq r\leq1.
        \end{cases}
    \end{equation*}
    Define the radial map $f_{M,\varepsilon}:\mathbb{B}^n\to\mathbb{B}^n$ by
    \begin{equation*}
        f_{M,\varepsilon}(\bfx)=
        \begin{cases}
            g_{M,\varepsilon}(\abs{\bfx})\dfrac{\bfx}{\abs{\bfx}},&\bfx\neq0,\\[4pt]
            0,&\bfx=0.
        \end{cases}
    \end{equation*}
    This is a smooth orientation-preserving self-diffeomorphism that equals the identity for $\abs{\bfx}\geq\varepsilon/2$. Its radial and tangential singular values are $g_{M,\varepsilon}'(r)$ and $g_{M,\varepsilon}(r)/r$, respectively, both positive and at most $M$. Since its global volume ratio is one, we obtain $E_{\ARAP}(f_{M,\varepsilon})\leq nM^2\mu(\mathbb{B}^n)\left({\varepsilon}/{2}\right)^n$.
    Meanwhile, $\nabla f_{M,\varepsilon}=MI$ for $\abs{\bfx}\leq\varepsilon/(4M)$, where the Jacobian determinant is $M^n$ and the volume-preserving energy density is $(M^n-1)^2$. Taking $\varepsilon=M^{-3/n}$ and letting $M\to\infty$ therefore yields vanishing integrated \ARAP\ distortion with unbounded local volume distortion. These observations motivate minimizing the joint angle- and volume-preserving energy, rather than the \ARAP\ energy directly, to better preserve the geometric properties of the parameterization.
\end{remark}

The preceding estimates concern integrated distortion and do not
by themselves enforce local nondegeneracy. To penalize local
collapse, we introduce the logarithmic \textbf{barrier energy}
\begin{equation}\label{eq:barrier_energy}
    E_B(f):=-\int_{\calM}
    \log\left(\Det(\nabla_{\calM}f)\right)\dif\mu.
\end{equation}
Its integrand tends to $+\infty$ as the local volume ratio
approaches zero. We therefore consider the penalized variational
problem
\begin{equation}\label{eq:barrier_optimization}
    \min_{f:\calM\to\mathcal N}
    E_\beta(f)+\eta E_B(f),
\end{equation}
where $\eta>0$ controls the strength of the barrier term, and $f$
satisfies the prescribed boundary and global bijectivity conditions.

\section{Discrete \boldmath$n$-manifold and energies}\label{sec:DiscreteEnergy}

In this section, we introduce discrete counterparts of the metric-preserving energies developed in the previous section for discrete $n$-manifolds and derive their differentials with respect to the mapping $f$. In the discrete setting, the map is not assumed a priori to be globally bijective
or {locally injective}, since folded simplices may occur during the numerical
optimization. Therefore, the simplex volumes used in the energy functionals below are understood as unsigned intrinsic volumes, while the exclusion of foldings is treated separately through the admissible class and boundary conditions. We begin by recalling the classical definitions of simplices, simplicial complexes, and discrete $n$-manifolds.

\begin{definition}\,
    \begin{enumerate}
        \item Let $\bfv_0,\dots,\bfv_k\in\mathbb{R}^n$ ($k\leq n$) be affinely independent. The $k$-simplex spanned by these vertices is $\sigma=[\bfv_0,\dots,\bfv_k]:=\left\{\sum_{i=0}^k\alpha_i\bfv_i:\alpha_i\geq0,\ \sum_{i=0}^k\alpha_i=1\right\}$.
        For $0\leq\ell<k$, its set of $\ell$-dimensional faces is $\Sigma_\ell(\sigma):=\left\{\tau\subseteq\sigma:\tau\text{ is spanned by $\ell+1$ vertices of }\sigma\right\}$.
        \item A simplicial complex $S$ is a collection of simplices that is closed under taking faces and such that the intersection of any two simplices is either empty or a common face. It is $n$-dimensional if all its maximal simplices are $n$-simplices.
        \item For a vertex $v\in S$, its link in $S$ is the subcomplex consisting of all simplices $\tau\in S$ such that $v\notin\tau$ and $\tau\cup\{v\}\in S$.

        \item A finite $n$-dimensional simplicial complex $\calM$ is called a discrete $n$-manifold, possibly with boundary, if its geometric realization is an $n$-dimensional piecewise-linear manifold, possibly with boundary. For $0\leq k\leq n$, the set of all $k$-simplices of $\calM$ is denoted by $\Sigma_k(\calM)$. In particular, $\Sigma_0(\calM)$ is the vertex set of $\calM$.
    \end{enumerate}
\end{definition}

Throughout this paper, $\mu$ denotes the intrinsic $n$-dimensional measure of $\calM$. For any $k$-simplex $\sigma$, the symbols $\abs{\sigma}$ and $\abs{f(\sigma)}$ denote the intrinsic $k$-dimensional Euclidean volumes of $\sigma$ and its image, respectively; the dimension $k$ is understood from context.

Let $N_v=\#(\Sigma_0(\calM))$ be the number of vertices in $\calM$. A piecewise affine map $f:\calM\to\re^d$ is uniquely determined by its values at these vertices, where $d\geq n$ is the dimension of the ambient Euclidean space containing the target manifold. We collect the vertex images in the matrix
\begin{equation*}
    \bff=\begin{bmatrix}
        f(\bfv_1)^\top\\
        \vdots\\
        f(\bfv_{N_v})^\top
    \end{bmatrix}
    \equiv\begin{bmatrix}
        \bff_1\\
        \vdots\\
        \bff_{N_v}
    \end{bmatrix}\in\re^{N_v\times d},
\end{equation*}
where $\bff_i=f(\bfv_i)^{\top}\in\re^{1\times d}$ denotes the image coordinates of the vertex $\bfv_i$. More precisely, on each $n$-simplex $\sigma=[\bfv_{i_0},\bfv_{i_1},\dots,\bfv_{i_n}]\in\Sigma_n(\calM)$, the map $f$ is given by the affine interpolation of its values at the vertices. Thus, for any point $\bfv\in\sigma$, $f(\bfv)=\sum_{j=0}^n\lambda_j(\bfv)\bff_{i_j}$ where $\lambda_j(\bfv)={\abs{[\bfv_{i_0},\dots,\bfv_{i_{j-1}},\bfv,\bfv_{i_{j+1}},\dots,\bfv_{i_n}]}}/{\abs{[\bfv_{i_0},\bfv_{i_1},\dots,\bfv_{i_n}]}}$ is the $j$th barycentric coordinate of $\bfv$ with respect to $\sigma$.

Unlike the continuous setting of the previous section, the discrete map $f$
is not assumed to be globally bijective or {locally injective} during the
numerical optimization. Let $\calM$ be oriented, and fix an orientation on the target $n$-manifold. For each nondegenerate oriented $n$-simplex $\tau\in\Sigma_n(\calM)$, let $\sgn f(\tau)\in\{-1,1\}$ denote the orientation sign of the affine restriction $f|_\tau$. We define the signed intrinsic $n$-dimensional volume of the image simplex by
\begin{equation}\label{eq:signed_simplex_volume}
    \svol(f(\tau)):=\sgn f(\tau)\cdot\abs{f(\tau)},
\end{equation}
where $\abs{f(\tau)}$ is the unsigned intrinsic $n$-dimensional Euclidean volume of $f(\tau)\subseteq\re^d$. Thus, $\svol(f(\tau))$ is positive on an orientation-preserving simplex and negative on an orientation-reversing simplex.
The signed local volume ratio on $\tau$ is then defined by
\begin{equation}\label{eq:signed_local_volume_ratio}
    J_\tau(f):=\frac{\svol(f(\tau))}{\abs{\tau}}.
\end{equation}
Throughout this section, the notation $\abs{\sigma}$ always denotes the unsigned intrinsic volume of a simplex $\sigma$. In particular, the volumes of all lower-dimensional faces are unsigned.

The following lemma follows directly from the \(n\)-dimensional
cotangent formula for the gradient of a simplex volume
\cite[Eq.~(3.17)]{ZHTL25}; see also \cite{Crane19}.

\begin{lemma}[\cite{Crane19,ZHTL25}]
\label{lem:gradient-signed-local-volume-ratio}
Let $\tau=[\bfv_{i_0},\ldots,\bfv_{i_n}]\in\Sigma_n(\calM)$ be an oriented \(n\)-simplex, and suppose that \(f(\tau)\) is nondegenerate. For \(j=0,\ldots,n\),
\begin{equation*}
  \nabla_{\bff_{i_j}}J_\tau(f)
  =
  \frac{\sgn f(\tau)}
       {n(n-1)\abs{\tau}}
  \sum_{\ell\neq j}\abs{f(\tau_{\widehat{i_j}\widehat{i_\ell}})}
  \cot\theta_{i_ji_\ell}^{\tau}(f)
  \bigl(\bff_{i_j}-\bff_{i_\ell}\bigr).
\end{equation*}
Here, \(\tau_{\widehat{i_j}\widehat{i_\ell}}\) denotes the
\((n-2)\)-face obtained by removing \(\bfv_{i_j}\) and
\(\bfv_{i_\ell}\) from \(\tau\), and
\(\theta_{i_ji_\ell}^{\tau}(f)\) is the corresponding dihedral angle
in \(f(\tau)\). For \(n=2\), the volume of a \(0\)-simplex is
understood to be one.
\end{lemma}

Since a discrete map may fold or overlap globally, the sum of the signed volumes of its image simplices is an algebraic, multiplicity-counted image volume. We denote it by
\begin{equation}\label{eq:discrete_image_volume}
    \SVol(f):=\sum_{\tau\in\Sigma_n(\calM)}\svol(f(\tau)).
\end{equation}
Accordingly, the discrete global volume ratio is defined by
\begin{equation}\label{eq:discrete_global_volume_ratio}
    \rmR_f:=\frac{\SVol(f)}{\mu(\calM)},\qquad \mu(\calM)=\sum_{\tau\in\Sigma_n(\calM)}\abs{\tau}.
\end{equation}
When $f$ is bijective and orientation-preserving, $\SVol(f)$ coincides with the intrinsic volume of the target image. In general, overlapping image simplices are counted with multiplicity, while orientation-reversing image simplices contribute negatively.

In light of the continuous metric-preserving energies
\cref{eq:Dirichlet_Energy,eq:Conformal_Energy,eq:Area_energy,eq:EVcontinuous,eq:EVScontinuous,eq:EVIDcontinuous},
we define the corresponding discrete energies as
\begin{equation*}
\left\{\begin{aligned}
    &E_D(f)
    =\sum_{\tau\in\Sigma_n(\calM)}\frac{\norm{\nabla_{\calM}f|_{\tau}}_{\rmF}^2}{n}\cdot\abs{\tau},\\
    &E_A(f)
    =\sum_{\tau\in\Sigma_n(\calM)}
    \abs{J_{\tau}(f)}^{2\over n}\cdot\abs{\tau},\\
    &E_C(f)=E_D(f)-E_A(f),
\end{aligned}\right.\qquad
\left\{\begin{aligned}
    &E_{V_s}(f)
    =\sum_{\tau\in\Sigma_n(\calM)}\left(J_{\tau}(f)\right)^2\cdot\abs{\tau},\\
    &E_{V_r}(f)={\SVol(f)^2\over\mu(\calM)},\\
    &E_V(f)=E_{V_s}(f)-E_{V_r}(f)
\end{aligned}\right.
\end{equation*}
for Dirichlet, area, conformal, volumetric stretch, volumetric ratio, and volume-preserving energies, respectively.
We also introduce the discrete joint angle- and volume-preserving energy $E_{\beta}(f)$ as in \cref{eq:Joint_Energy} and define the discrete logarithmic barrier energy, inspired by \cref{eq:barrier_energy}, by
\begin{equation}\label{eq:discrete_barrier_energy}
    E_B(f)
    =-\sum_{\tau\in\Sigma_n(\calM)}
    \log J_{\tau}(f)\cdot\abs{\tau},
    \qquad f\in\mathcal A_+,
\end{equation}
where $\mathcal{A}_{+}:=\{f:J_\tau(f)>0,\ \forall\,\tau\in\Sigma_n(\calM)\}$. The condition $f\in\mathcal{A}_{+}$ {enforces local injectivity on every simplex},
which excludes simplex inversion and degeneration. The additional topological conditions
needed for global bijectivity are stated in \Cref{thm:global-bijectivity-certificate}.

To reveal the common algebraic structure of the discrete metric-preserving energies and their gradients, we express their constituent terms through weighted Laplacian matrices. For symmetric edge weights $w=\{w_{ij}\}_{i\neq j}$, define
\begin{equation*}
  [\mathcal L(w)]_{ij}=
  \begin{cases}
    -w_{ij},
      &i\neq j\ \text{and}\ [\bfv_i,\bfv_j]\in\Sigma_1(\calM),\\
    \displaystyle\sum_{\ell\neq i}w_{i\ell},
      &i=j,\\
    0,
      &\text{otherwise}.
  \end{cases}
\end{equation*}
For an edge $[\bfv_i,\bfv_j]\subset\tau$, introduce the source and image
cotangent kernels
\begin{equation*}
  \omega_{D,ij}^{\tau}
  :=\frac{\abs{\tau_{\hat i\hat j}}\cot\theta_{ij}^{\tau}}
  {n(n-1)},
  \qquad
  \omega_{ij}^{\tau}(f)
  :=\frac{\abs{f(\tau_{\hat i\hat j})}
  \cot\theta_{ij}^{\tau}(f)}{n(n-1)}.
\end{equation*}
Here $\tau_{\hat i\hat j}$ is the $(n-2)$-face obtained by removing
$\bfv_i$ and $\bfv_j$. The angles $\theta_{ij}^{\tau}$ and
$\theta_{ij}^{\tau}(f)$ are the corresponding intrinsic dihedral angles in
$\tau$ and $f(\tau)$, respectively. When $n=2$, a $0$-simplex is assigned
unit $0$-dimensional volume. Set $L_D=\mathcal L(w_D)$ and $L_Q(f)=\mathcal L(w_Q(f))$ for $Q\in\{A,V_s,V_r,B\}$, where $w_{D,ij}=\sum_{\tau\supset[\bfv_i,\bfv_j]}\omega_{D,ij}^{\tau}$, $w_{Q,ij}(f)=\sum_{\tau\supset[\bfv_i,\bfv_j]}\gamma_Q^{\tau}(f)\omega_{ij}^{\tau}(f)$.
The scalar multipliers $\gamma_Q^{\tau}$ are listed in the second column of \cref{tab:sec3-energy-gradient-summary}.

\begin{table}[t]
  \centering
  \caption{Simplexwise multipliers and Laplacian representations of the discrete energies.}
  \label{tab:sec3-energy-gradient-summary}
  \scriptsize
  \renewcommand{\arraystretch}{1.45}
  \begin{tabular}{cccc}
    \toprule
    Energy $E_Q$ & $\gamma_{Q}^{\tau}(f)$ & Discrete form & Gradient w.r.t. $\bff$\\
    \midrule
    $E_D$ & --
      & $\displaystyle
        \frac1n\tr(\bff^\top L_D\bff)$
      & $\displaystyle
        \frac2nL_D\bff$\\
    $E_A$ & $\displaystyle\frac{\abs{f(\tau)}^{\frac{2}{n}-1}}{\abs{\tau}^{\frac{2}{n}-1}}$
      & $\displaystyle
        \sum_{\tau\in\Sigma_n(\calM)}
        \frac{\abs{f(\tau)}^{2/n}}
        {\abs{\tau}^{2/n-1}}
        =\frac1n\tr(\bff^\top L_A(f)\bff)$
      & $\displaystyle
        \frac2nL_A(f)\bff$\\
    $E_{V_s}$ & $\displaystyle\frac{\abs{f(\tau)}}{\abs{\tau}}$
      & $\displaystyle
        \sum_{\tau\in\Sigma_n(\calM)}
        \frac{\abs{f(\tau)}^2}{\abs{\tau}}
        =\frac1n\tr(\bff^\top L_{V_s}(f)\bff)$
      & $\displaystyle
        2L_{V_s}(f)\bff$\\
    $E_{V_r}$ & $\displaystyle\rmR_f\sgn f(\tau)$
      & $\displaystyle
        \frac1n\tr(\bff^\top L_{V_r}(f)\bff)$
      & $\displaystyle
        2L_{V_r}(f)\bff$\\
    $E_B$ & $\displaystyle\frac{\abs{\tau}}{\svol(f(\tau))}$
      & $\displaystyle
        -\sum_{\tau\in\Sigma_n(\calM)}
        \abs{\tau}\log J_\tau(f)$
      & $\displaystyle
        -L_B(f)\bff$\\
    \bottomrule
  \end{tabular}
\end{table}

\begin{theorem}\label{thm:discrete_energies_gradients}
    Let $\calM$ be an oriented discrete $n$-manifold with $n\geq2$, and let $f:\calM\to\re^d$, $d\geq n$, be piecewise affine on every $n$-simplex. Then the constituent energies $E_D$, $E_A$, $E_{V_s}$, and $E_{V_r}$, together with their gradients with respect to $\bff$, admit the Laplacian representations summarized in \cref{tab:sec3-energy-gradient-summary}. For $f\in\mathcal A_+$, the corresponding representation for $E_B$ is given in the last row of the same table. Define
    \begin{equation*}
        L_C(f):=L_D-L_A(f),\qquad L_V(f):=L_{V_s}(f)-L_{V_r}(f).
    \end{equation*}
    Then the conformal and volume-preserving energies defined above satisfy
    \begin{align*}
        &E_C(f)=\frac1n\tr\left(\bff^\top L_C(f)\bff\right),\qquad\nabla_{\bff}E_C(f)=\frac2nL_C(f)\bff,\\
        &E_V(f)=\frac1n\tr\left(\bff^\top L_V(f)\bff\right),\qquad\nabla_{\bff}E_V(f)=2L_V(f)\bff.
    \end{align*}
\end{theorem}
\begin{proof}
    The $n$-dimensional cotangent identity and its differential 
    \cite{Crane19,ZHTL25} give 
    \begin{equation}\label{eq:proof-simplex-cotangent-identity} 
    \abs{f(\tau)} 
    =\frac1n\sum_{[\bfv_i,\bfv_j]\subset\tau} 
    \omega_{ij}^{\tau}(f)\norm{\bff_i-\bff_j}^{2}, 
    \qquad 
    \nabla_{\bff_i}\abs{f(\tau)} 
    =\sum_{\bfv_j\in\tau,\,j\ne i} 
    \omega_{ij}^{\tau}(f)(\bff_i-\bff_j). 
    \end{equation} 
    The second identity also follows from 
    \cref{lem:gradient-signed-local-volume-ratio}, since the orientation sign is 
    locally constant on the set of nondegenerate simplices. Applying 
    \cref{eq:proof-simplex-cotangent-identity} and the chain rule to 
    ${\abs{f(\tau)}^{2/n}}/{\abs{\tau}^{2/n-1}}$ gives the multiplier 
    $\gamma_A^\tau$, and hence the stated energy and gradient formulas for 
    $E_A$. The formulas for $E_D$ and $E_{V_s}$ follow directly from 
    \cite{ZHTL25}. Subtraction therefore 
    proves the assertions for $E_C$. 
     
    It remains to account for the volume and barrier terms. Since, by \cref{eq:signed_simplex_volume}, $\svol(f(\tau))
    =\operatorname{sgn} f(\tau)\abs{f(\tau)}$, 
    \cref{eq:proof-simplex-cotangent-identity} yields 
    \[\nabla_{\bff}\SVol(f)=\mathcal L\left( \left\{\sum_{\tau\supset[\bfv_i,\bfv_j]} \operatorname{sgn} f(\tau)\omega_{ij}^{\tau}(f)\right\}_{ij}\right)\bff.\] 
    
    Thus, using \cref{eq:discrete_image_volume,eq:discrete_global_volume_ratio},
    \(\rmR_f=\SVol(f)/\mu(\calM)\) and
    \(E_{V_r}(f)=\SVol(f)^2/\mu(\calM)\),
    \[\nabla_{\bff}E_{V_r}(f)=2\rmR_f\nabla_{\bff}\SVol(f)=2L_{V_r}(f)\bff.\]
    The cotangent identity also gives $\frac1n\tr\!\left(\bff^\top L_{V_r}(f)\bff\right)=\rmR_f\SVol(f)=E_{V_r}(f)$.
    Moreover, 
    \[\sum_\tau\abs{\tau}\bigl(J_\tau(f)-\rmR_f\bigr)^2=\sum_\tau\frac{\abs{f(\tau)}^2}{\abs{\tau}}-\frac{\SVol(f)^2}{\mu(\calM)}=E_{V_s}(f)-E_{V_r}(f),\] 
    which proves all statements for \(E_V\). By \cref{lem:gradient-signed-local-volume-ratio}, 
    \[\begin{aligned} 
      \nabla_{\bff_i}E_B(f) 
      =-\sum_{\tau\ni\bfv_i} 
        \frac{\abs{\tau}}{\abs{f(\tau)}} 
        \sum_{\bfv_j\in\tau,\,j\ne i} 
        \omega_{ij}^{\tau}(f)(\bff_i-\bff_j)
      =-[L_B(f)\bff]_i. 
    \end{aligned}\]
    This is exactly the multiplier 
    $\gamma_B^\tau={\abs{\tau}}/{\svol(f(\tau))}$ on \(\mathcal A_+\). 
    Finally, $-\abs{\tau}\log J_\tau(f)\to+\infty$ as $J_\tau(f)\downarrow0$. 
    Hence, \(E_B\) diverges as the boundary of \(\mathcal A_+\) is approached from 
    the orientation-preserving side.
\end{proof}

The only distinction among the Laplacian matrices associated with the constituent energies is the simplexwise scalar multiplier listed in \cref{tab:sec3-energy-gradient-summary}. With $E_B(f)=+\infty$ outside $\mathcal A_+$, the barrier tends to $+\infty$ whenever an image simplex approaches degeneracy from the orientation-preserving side.

\section{A unified numerical framework for bijective metric-preserving parameterizations}
\label{sec:alg}

Building on the energies and discrete gradient formulas developed in
\Cref{sec:ContinuousEnergy,sec:DiscreteEnergy}, this section presents a unified three-stage numerical framework for computing bijective joint angle- and volume-preserving parameterizations. Stage~I constructs a low-energy map satisfying the target geometry, Stage~II restores {local injectivity} by ADMM with an inexact PNCG solution of the mapping subproblem, and Stage~III solves the discrete counterpart of \cref{eq:barrier_optimization} by minimizing
\begin{equation}\label{eq:discrete_joint_energy}
    E_{\beta,\eta}(f):=E_{\beta}(f)+\eta E_B(f),\qquad\beta\in[0,1],\,\eta>0,
\end{equation}
which is a weighted sum of the joint angle- and volume-preserving energy $E_{\beta}(f)$ in \cref{eq:Joint_Energy} and the logarithmic barrier energy $E_B(f)$ in \cref{eq:discrete_barrier_energy}. Since $E_B$ is not defined when the initial map contains degenerate or inverted simplices, Stage~II precedes the repairing procedure whenever folding is present.

The framework is formulated for parameterizations between general $n$-dimensional manifolds. In what follows, we develop spherical and ball parameterizations in detail, with the corresponding geometric settings summarized in \cref{tab:sec4-problem-instances}. In these two settings, $\mathtt B$ denotes the set of vertices constrained to the unit sphere. Apart from the initialization in Stage~I, both problems use the same optimization framework in Stages~II and~III; their differences are encoded by $(\calM,n,d,\mathtt B)$ and the associated sphere-constrained vertices.
\begin{table}
  \centering
  \caption{Two instances of the unified numerical framework.}
  \label{tab:sec4-problem-instances}
  \footnotesize
  \begin{tabular}{ccccc}
    \toprule
    Problem & Source $\calM$ & $d$ & Target $\calN\subset\re^d$
      & $\mathtt B$\\
    \midrule
    Spherical parameterization
      & $\calM\simeq\bbS^n$ & $n+1$ & $\bbS^n$
      & $\Sigma_0(\calM)$\\
    Ball parameterization
      & $\calM\simeq\bbB^n$ & $n$ & $\bbB^n$
      & $\Sigma_0(\partial\calM)$\\
    \bottomrule
  \end{tabular}
\end{table}

{
\begin{theorem}
\label{thm:global-bijectivity-certificate}
Let $\calM$ be a connected, oriented finite simplicial $n$-manifold.
\begin{enumerate}[label=\textnormal{(\roman*)}]
\item Suppose $\partial\calM=\varnothing$ and let
$f:\calM\to\re^{n+1}$ be piecewise affine with $0\notin f(\calM)$.
Define $\widehat f=f/\norm{f}:\calM\to\bbS^n$.
If $\widehat f$ is locally injective and $\deg(\widehat f)=1$,
then $f:\calM\to f(\calM)$ and $\widehat f:\calM\to\bbS^n$
are homeomorphisms.

\item Suppose $\calM$ is homeomorphic to the closed $n$-ball.
If $f:\calM\to\mathbb R^n$ is piecewise affine,
$J_\tau(f)>0$ for every $\tau\in\Sigma_n(\calM)$, and
$f|_{\partial\calM}$ is a homeomorphism onto the boundary of a bounded
polyhedral domain $\Omega$, then $f$ is a global homeomorphism onto
$\overline{\Omega}$.
\end{enumerate}
\end{theorem}

\begin{proof}
In (i), classical covering-space and degree theory imply that
$\widehat f$ is a one-sheeted covering of $\bbS^n$, hence a
homeomorphism~\cite{Hatcher2002}. Consequently,
$f:\calM\to f(\calM)$ is also a homeomorphism.
Part~(ii) follows from the simplicial global inversion theorem of
Lipman~\cite{Lipman2014}.
\end{proof}
}

\Cref{thm:global-bijectivity-certificate} shows that, to obtain a globally
bijective spherical parameterization, it suffices to ensure that
$\widehat f$ has degree one and is locally injective. For ball
parameterizations, it suffices to ensure boundary bijectivity and positive
simplex Jacobians. The analysis below establishes the degree-one condition
under an approximation assumption, while Stages~II and~III remove foldings
and maintain positive simplex Jacobians.

{
In (i), $\widehat f$ must be injective on a neighborhood of each point,
including points on simplex interfaces. A sufficient discrete condition,
which will be used for the a posteriori bijectivity checks in
\Cref{sec:numerical-experiments}, is that $\widehat f$ be an embedding on the closed star of every vertex
~\cite[Lemma~15]{DyerVegterWintraecken2015}.
}

\subsection{Stage I: Initialization}
\label{subsec:stage-I-initialization}

We give explicit initialization procedures for the targets $\bbS^n$ and $\bbB^n$. Before initialization, the source mesh is uniformly scaled so that its total intrinsic volume agrees with that of the target. If the source volume before scaling is $\mu(\calM)$, we set $c_{\calM}=\left({\nu(\calN)}/{\mu(\calM)}\right)^{1\over n}$, $\mu(c_{\calM}\calM)=\nu(\calN)$. The scaled source mesh $c_{\calM}\calM$ is again denoted by $\calM$ for simplicity. With this normalization, the global volume ratio in \cref{eq:discrete_global_volume_ratio} becomes $\rmR_f=\SVol(f)/\nu(\calN)\approx1$ for an orientation-preserving bijective parameterization whose image approximates the target.

\paragraph{North--south initialization for spherical parameterization}
\label{subsec:spherical-parameterization}

Let $\calM$ be an oriented simplicial $n$-manifold homeomorphic to $\bbS^n$.
The target is $\calN=\bbS^n\subset\re^{n+1}$, and we seek a piecewise-affine
map $f:\calM\to\re^{n+1}$ with $f(\bfv)\in\bbS^n$ for every
$\bfv\in\Sigma_0(\calM)$.

In light of \Cref{thm:global-bijectivity-certificate}, we seek an initial
map $f^{(0)}$ whose radial normalization has degree one. The following
theorem establishes this property for the continuous generalized Dirac map.
Its discrete counterpart extends the surface construction of
Haker et al.~\cite{SHSA00} and was used by Tan et al.~\cite{ZHTL25}
to initialize $n$-dimensional spherical parameterizations.
We take this discrete Dirac map as $f^{(0)}$ and refine it through the
north--south iteration described below.

\begin{theorem}
\label{thm:continuous-dirac-degree}
Let $\calM$ be a closed, connected, oriented smooth $n$-manifold homeomorphic
to $\bbS^n$, $n\geq2$, and let $\mathbf{p}\in \calM$. In positively oriented geodesic
normal coordinates $u=(u_1,\ldots,u_n)$ at $\mathbf{p}$, let $h=(h_1,\ldots,h_n)$ be the zero-mean distributional solution of $-\Delta_{\calM} h_i=\partial_{p_i}\delta_{\mathbf{p}}$, ($i=1,\ldots,n$), where
$\langle\partial_{p_i}\delta_{\mathbf{p}},\varphi\rangle
=\partial_{u_i}\varphi(\mathbf{p})$.
Let $\Pi:\bbS^n\setminus\{N\}\to\mathbb R^n$ be stereographic projection
and define $g=\Pi^{-1}\circ h$ on $\calM\setminus\{\mathbf{p}\}$ with $g(\mathbf{p})=N$.
With compatible orientations, it holds that $\deg(g)=1$.
\end{theorem}
\begin{proof}
Let $G(\bfx,\mathbf{p})$ be the zero-mean Green function of $-\Delta_{\calM}$,
which satisfies, in the distributional sense with respect to $\bfx$, $-\Delta_{\calM}G(\bfx,\mathbf{p})=\delta_{\mathbf{p}}-{1}/{\mu(\calM)}$ and $\int_{\calM}G(\bfx,\mathbf{p})\dif\mu(\bfx)=0$.
Then $h_i(\bfx)=\partial_{p_i}G(\bfx,\mathbf{p})$, and the local Green-function
expansion~\cite[pp.~101--114]{Aubin1998} gives $h(u)=c_n{u}/{|u|^n}+o(|u|^{1-n})$ for some $c_n>0$.
Hence $|h(u)|\to\infty$ as $u\to0$, so $g$ extends continuously to
$g(\mathbf{p})=N$ and $g^{-1}(N)=\{\mathbf{p}\}$. In an inverted local coordinate at $N$, ${h(u)}/{|h(u)|^2}=c_n^{-1}|u|^{n-2}u+o(|u|^{n-1})$.
After radial normalization, its restriction to a sufficiently small
surrounding sphere is homotopic to the identity on $\bbS^{n-1}$. So the local degree at $\mathbf{p}$ is $+1$. Since $\mathbf{p}$ is the unique preimage of $N$, $\deg(g)=1$.
\end{proof}

\begin{lemma}
\label{lem:degree-preservation}
Let $\calM$ be a closed, connected, oriented finite simplicial $n$-manifold. Assume $g:\calM\to\bbS^n$ is continuous, and $f:\calM\to\re^{n+1}$ is piecewise affine with $f(\bfv)\in\bbS^n$ for every $\bfv\in\Sigma_0(\calM)$. Set $\delta=\max\limits_{\bfv\in\Sigma_0(\calM)}\norm{f(\bfv)-g(\bfv)}$, $\omega_g=\max\limits_{\tau\in\Sigma_n(\calM)}\max\limits_{\bfv\in\Sigma_0(\tau)}\sup\limits_{\bfx\in\tau}\norm{g(\bfx)-g(\bfv)}$. Then $f$ never vanishes, $\widehat f:=f/\norm{f}:\calM\to\bbS^n$ is homotopic to $g$, and $\deg(\widehat f)=\deg(g)$, provided one of the following conditions holds: (i) $\delta^2/2+\omega_g<1$; (ii) For some $0<\varepsilon<\sqrt{2}$, there exists $h_0=h_0(g,\varepsilon)>0$, such that if every simplex of $\calM$ has Euclidean diameter less than $h_0(g,\varepsilon)$, and $\norm{f(\bfv)-g(\bfv)}<\sqrt{2}-\varepsilon$ for every $\bfv\in\Sigma_0(\calM)$.
\end{lemma}
\begin{proof}
For (i), let $\bfx\in\tau$ and $\bfv\in\Sigma_0(\tau)$. The unit lengths of
$f(\bfv)$ and $g(\bfv)$ give
\begin{equation*}
    f(\bfv)\cdot g(\bfx)=1-\tfrac12\norm{f(\bfv)-g(\bfv)}^2+f(\bfv)\cdot\bigl(g(\bfx)-g(\bfv)\bigr)\geq1-\delta^2/2-\omega_g>0.
\end{equation*}
Since $f(\bfx)$ is a convex combination of the vertex values,
$f(\bfx)\cdot g(\bfx)>0$, so $f(\bfx)\neq0$.
Moreover, $((1-t)g(\bfx)+tf(\bfx))\cdot g(\bfx)>0$ for every
$t\in[0,1]$. Thus $H(\bfx,t):=((1-t)g(\bfx)+tf(\bfx))/{\norm{(1-t)g(\bfx)+tf(\bfx)}}$
is a continuous homotopy from $g$ to $\widehat f$. The degree is
preserved.

For (ii), uniform continuity of $g$ on the compact source $\calM$ gives $h_0>0$
such that $\norm{g(\bfx)-g(\bfy)}<1-\tfrac12(\sqrt{2}-\varepsilon)^2$
whenever $\bfx,\bfy\in\calM$ and $\norm{\bfx-\bfy}<h_0$.
If every simplex has Euclidean diameter less than $h_0$, then
$\omega_g\leq1-\tfrac12(\sqrt{2}-\varepsilon)^2$.
Together with $\delta<\sqrt{2}-\varepsilon$, this gives
$\delta^2/2+\omega_g<1$, so (i) applies.
\end{proof}

{
For the spherical problem, let $g$ be the continuous generalized Dirac
map in \Cref{thm:continuous-dirac-degree}, regarded on the triangulated
source through an orientation-preserving homeomorphism.
For any map $f$ obtained during subsequent refinement, if its vertexwise
error $\delta$ relative to $g$ satisfies $\delta^2/2+\omega_g<1$, then
\Cref{lem:degree-preservation} shows that $\widehat f$ is homotopic to $g$.
Consequently, $\deg(\widehat f)=\deg(g)=1$.
This degree-one property is preserved in subsequent iterations if every
vertex displacement is at most $\sqrt{2}-\varepsilon$ for a fixed
$0<\varepsilon<\sqrt{2}$, provided the triangulation is sufficiently fine
relative to the radial normalization of each current iterate.}

By \Cref{sec:DiscreteEnergy}, the joint angle- and volume-preserving energy can be decomposed as
\begin{equation*}
    E_\beta(f)=\beta E_D(f)+(1-\beta)E_{V_s}(f)-\beta E_A(f)-(1-\beta)E_{V_r}(f).
\end{equation*}
To construct the Stage-I initialization, we temporarily omit the two negative correction terms $E_A$ and $E_{V_r}$ and use only the positive part $\widehat{E}(f):=\beta E_D(f)+(1-\beta)E_{V_s}(f)$. Given the current iterate $f^{(k)}$, we further freeze the map-dependent Laplacian associated with $E_{V_s}$ at $f^{(k)}$. The resulting iteration matrix is
\begin{equation}
    \widehat L^{(k)}\equiv\widehat{L}(f^{(k)}):=\frac{\beta}{n}L_D+(1-\beta)L_{V_s}(f^{(k)}).\label{eq:sec4-stage-I-matrix}
\end{equation}
Here $L_D$ is the fixed source finite element Laplacian, whereas $L_{V_s}(f^{(k)})$ is the map-dependent matrix associated with $E_{V_s}$ in \Cref{sec:DiscreteEnergy}.

This positive-term lagging strategy is motivated by existing spherical SEM/VSEM methods. On a well-shaped two-dimensional Delaunay mesh, the standard cotangent matrix has an M-matrix structure and the associated harmonic update has a convex-combination property that favors orientation preservation~\cite{MHTL19}. North--south alternating solutions of the Laplacian systems also have R-linear convergence under the corresponding assumptions~\cite{TMWH24}. Thus, \cref{eq:sec4-stage-I-matrix} retains the stability and efficiency of a sparse Laplacian system.

Moreover, once the image of $f$ approximately covers the fixed target
$\calN$, the two negative correction terms vary slowly. With the discrete
energies of \Cref{sec:DiscreteEnergy}, when $f$ is orientation preserving, has little overlap,
and approximately covers $\calN$, for $n=2$ we have
\begin{equation*}
  E_A(f)
  =\sum_{\tau\in\Sigma_2(\calM)}\abs{f(\tau)}
  \approx\nu(\calN),
  \qquad
  E_{V_r}(f)
  =\frac{\SVol(f)^2}{\mu(\calM)}
  \approx\frac{\nu(\calN)^2}{\mu(\calM)}.
\end{equation*}
Hence, near full coverage, $E_A$ and $E_{V_r}$ contribute only small variations and gradient corrections, which motivates using the positive terms to generate the initialization direction.

Let $\Pi:\bbS^n\setminus\{(0,\ldots,0,1)\}\to\re^n$ be the fixed stereographic projection used in the implementation, and let $\Pi^{-1}$ denote its inverse. At iteration $k$, we assemble $\widehat L^{(k)}$ from $f^{(k)}$ and keep it fixed. We project the current vertex images and apply radial inversion in place: $\bfy_i=\Pi(\bff_i^{(k)})$, $\bfy_i\leftarrow{\bfy_i}/{\norm{\bfy_i}^2}$ ($i=1,\ldots,N_v$).
For a prescribed chart radius $r>0$, set
$\mathtt I_k=\{i:\norm{\bfy_i}<r\}$ and
$\mathtt B_k=\{1,\ldots,N_v\}\setminus\mathtt I_k$.
We update the vertices in this bounded projection region, away from the projection pole, while the remaining vertices provide Dirichlet data. Keeping $\bfy_{\mathtt B_k}$ fixed, we solve $[\widehat L^{(k)}]_{\mathtt I_k\mathtt I_k}\bfy_{\mathtt I_k}=-[\widehat L^{(k)}]_{\mathtt I_k\mathtt B_k}\bfy_{\mathtt B_k}$.
The next iterate is given by $\bff_i^{(k+1)}=\Pi^{-1}(\bfy_i)$ for $i=1,\ldots,N_v$. The inverse projection already places each vertex image on $\bbS^n$; the implementation additionally normalizes these values to correct roundoff errors.

Since \cref{eq:sec4-stage-I-matrix} uses only the current coefficients of the
positive terms, the north--south iteration does not exactly minimize
$E_\beta$. It is used only to construct a spherical starting map for
Stages~II--III, and therefore only a moderate initialization accuracy is
required. After $K_{\mathrm I}$ iterations, the resulting map is denoted by
$f^{\mathrm I}$. Its vertex images lie on $\bbS^n$, but its image simplices
may still be degenerate or inverted.

{\begin{algorithm}
\caption{North--south iteration for spherical initialization}
\label{alg:sec4-north-south-initialization}
\begin{algorithmic}[1]
\REQUIRE The $n$-manifold $\calM$, $\beta\in[0,1]$, a chart radius $r>0$, and number of iterations $K_{\mathrm I}$
\ENSURE A piecewise-affine map $f^{\mathrm I}:\calM\to\re^{n+1}$
with $f^{\mathrm I}(\bfv)\in\bbS^n$ for every $\bfv\in\Sigma_0(\calM)$
\STATE Compute a Dirac initial map $\bff^{(0)}$
\FOR{$k=0,\ldots,K_{\mathrm I}-1$}
  \STATE Assemble
  $\widehat L^{(k)}=\frac{\beta}{n}L_D
  +(1-\beta)L_{V_s}(\bff^{(k)})$
  \STATE Set $\bfy_i=\Pi(\bff_i^{(k)})$ for $i=1,\ldots,N_v$
  \STATE Update $\bfy_i\leftarrow\bfy_i/\norm{\bfy_i}^2$
  for $i=1,\ldots,N_v$
  \STATE Set
  $\mathtt I_k=\{i:\norm{\bfy_i}<r\}$ and
  $\mathtt B_k=\{1,\ldots,N_v\}\setminus\mathtt I_k$
  \STATE Keeping $\bfy_{\mathtt B_k}$ fixed, solve
  $[\widehat L^{(k)}]_{\mathtt I_k\mathtt I_k}
  \bfy_{\mathtt I_k}
  =-[\widehat L^{(k)}]_{\mathtt I_k\mathtt B_k}
  \bfy_{\mathtt B_k}$
  \STATE Set $\bff_i^{(k+1)}=\Pi^{-1}(\bfy_i)$
  for $i=1,\ldots,N_v$
\ENDFOR
\STATE Set $f^{\mathrm I}=f^{(K_{\mathrm I})}$
\end{algorithmic}
\end{algorithm}}

\paragraph{Boundary map and harmonic extension for ball parameterization}
\label{subsec:ball-parameterization}

Let $\calM$ be an oriented simplicial $n$-manifold homeomorphic to $\bbB^n$, with $\partial\calM\simeq\bbS^{n-1}$. Then $\calN=\bbB^n\subset\re^n$, $d=n$, and we seek a piecewise-affine map $f:\calM\to\re^n$ satisfying $f(\calM)\subset\bbB^n$ and $f(\bfv)\in\bbS^{n-1}$ for every $\bfv\in\Sigma_0(\partial\calM)$.
We first compute a {locally injective} piecewise-affine boundary map
$f_\partial:\partial\calM\to\re^n$, with
$f_\partial(\bfv)\in\bbS^{n-1}$ for every $\bfv\in\Sigma_0(\partial\calM)$.
This is a spherical subproblem of intrinsic dimension $n-1$; applying the preceding
spherical initialization and Stages~II--III in that dimension yields $f_\partial$.
{This boundary construction prepares for the subsequent application of the
ball bijectivity criterion in \Cref{thm:global-bijectivity-certificate}.}
For $n=2$, we preserve the cyclic order of the boundary vertex images on
$\bbS^1$; for $n\geq3$, boundary bijectivity can be checked using the spherical
criterion of that theorem.
Let $\mathtt I$ and $\mathtt B$ denote the interior and boundary vertex sets, respectively. Set $\bff_{\mathtt B}^{\mathrm I}=\bff_\partial$ and initialize the interior vertices by the discrete harmonic extension $[L_D]_{\mathtt I\mathtt I}\bff_{\mathtt I}^{\mathrm I}=-[L_D]_{\mathtt I\mathtt B}\bff_{\mathtt B}^{\mathrm I}$. The resulting map $f^{\mathrm I}$ has boundary vertex images on $\bbS^{n-1}$ and low Dirichlet energy, but its image $n$-simplices are not necessarily {locally injective}.

\subsection{Stage II: Folding repair}
\label{subsec:stage-II-folding-repair}

The Stage-I map satisfies the sphere constraints at the designated vertices but may contain degenerate or inverted image simplices. Stage~II repairs these foldings while maintaining the sphere constraints and produces a strictly feasible starting map for Stage~III.

Fix $\varepsilon>0$ and define
$\mathcal C_{\varepsilon}
=\left\{z:z_\tau\geq\varepsilon,\ \tau\in\Sigma_n(\calM)\right\}$.
We consider the proximal ADMM framework for the following
optimization problem
\begin{equation}
\begin{aligned}
    \min_{\bff,z}\quad&
    I_{\mathcal C_{\varepsilon}}(z)
    =\begin{cases}
        0,&z\in\mathcal C_{\varepsilon},\\
        +\infty,&z\notin\mathcal C_{\varepsilon},
    \end{cases}
    \qquad\mathrm{s.t.}\quad
    \begin{aligned}
        J_\tau(f)&=z_\tau,\ \tau\in\Sigma_n(\calM),\\
        \norm{\bff_i}&=1,\ i\in\mathtt B.
    \end{aligned}
\end{aligned}
\label{eq:sec4-repair-model}
\end{equation}
The constraint $z_\tau=J_\tau(f)$, together with
$z\in\mathcal C_{\varepsilon}$, enforces
$J_\tau(f)\geq\varepsilon$ for every $\tau\in\Sigma_n(\calM)$.
Consequently, any pair $(\bff,z)$ satisfying these constraints and the sphere
constraints is a global minimizer of \cref{eq:sec4-repair-model} with zero
objective value. Since such minimizers are generally nonunique, we initialize
Stage~II with $\bff^{(0)}=\bff^I$ and employ proximal regularization to
discourage large changes during the folding-repair process. The proximal term $\norm{\bff-\bff^{(k)}}_{\mathrm{F}}^2$ controls the
displacement of an ADMM iteration and stabilizes the mapping subproblem.
For spherical maps, limiting vertex displacements is also motivated by
\Cref{lem:degree-preservation}: any iterate satisfying the lemma's approximation
condition relative to the continuous Dirac map has a radial normalization
of degree one.
Using the proximal ADMM framework, one inexact iteration is
\begin{align}
    z_\tau^{(k+1)}
    &=\underset{z_\tau\ge\varepsilon}{\arg\min}\ 
    \frac{1}{2}\bigl(J_\tau(f^{(k)})-z_\tau+u_\tau^{(k)}\bigr)^2
    =
    \max\{\varepsilon,
    J_\tau(f^{(k)})+u_\tau^{(k)}\},
    \notag\\
  \bff^{(k+1)}
  &\approx\underset{\bff\in\calX_{\mathtt B}}{\arg\min}\ \left\{\phi_k(f):=\sum_{\tau\in\Sigma_n(\calM)}\frac{\abs{\tau}}{2}
  \bigl(J_\tau(f)-z_\tau^{(k+1)}+u_\tau^{(k)}\bigr)^2
  +\frac{\kappa}{2}\norm{\bff-\bff^{(k)}}_{\mathrm{F}}^2\right\},
  \label{eq:sec4-f-update}\\
  u_\tau^{(k+1)}
  &=u_\tau^{(k)}+J_\tau(f^{(k+1)})-z_\tau^{(k+1)},
  \notag
\end{align}
where
$\calX_{\mathtt B}
=\left\{\bff\in\re^{N_v\times d}:
\norm{\bff_i}=1,\ i\in\mathtt B\right\}$, $u_\tau$ ($\tau\in\Sigma_n(\calM)$) are scaled dual variables, and $\kappa>0$ is the proximal regularization parameter.

The mapping subproblem is a constrained optimization problem on
$\calX_{\mathtt B}$. Thus, all vertices are updated tangentially in the spherical problem, whereas
the ball problem uses Euclidean updates for interior vertices and tangential
updates for boundary vertices.

Let $\varrho_\tau^{(k)}(f)=J_\tau(f)-z_\tau^{(k+1)}+u_\tau^{(k)}$. The gradient formula in \cref{lem:gradient-signed-local-volume-ratio} gives $\nabla_{\bff_i}\phi_k(f)={\sum_{\bfv_i\in\tau\in\Sigma_n(\calM)}\abs{\tau}\varrho_\tau^{(k)}(f)\nabla_{\bff_i}J_\tau(f)}+\kappa(\bff_i-\bff_i^{(k)})$, and the constrained gradient used by the inner solver is $\operatorname{grad}\phi_k(f)=\mathcal P_{\bff}\!\left(\nabla_{\bff}\phi_k(f)\right)$. Here, \(\mathcal P_{\bff}\) denotes the orthogonal projection onto the tangent space of the constraint manifold at \(\bff\). Specifically, it removes the radial component at sphere-constrained vertices and acts as the identity on unconstrained vertices.
We solve \cref{eq:sec4-f-update} inexactly by a restarted preconditioned
Riemannian nonlinear conjugate-gradient method~\cite{AMS08,NW06}. We refer to this proximal ADMM scheme with an inexact PNCG solution of the
mapping subproblem as the ADMM--PNCG folding-repair stage. The preconditioner is the diagonal Gauss--Newton approximation of the squared residual term in \cref{eq:sec4-f-update}. Trial maps preserve the sphere constraints at constrained vertices and use ordinary Euclidean updates elsewhere, with step sizes determined by an Armijo backtracking line search. Stage~II is successful only after the direct verification $\min_{\tau\in\Sigma_n(\calM)}J_\tau(f^{\mathrm{II}})>0$.

\begin{remark}
    In practical implementations, it is unnecessary to update all vertices during Stage~II,
    since foldings are typically confined to a small portion of the mesh.
    We therefore update only the vertices of simplices with nonpositive or small
    positive Jacobians and a small surrounding neighborhood,
    while keeping the remaining vertices fixed. The repair region is adjusted
    as the iterations proceed. The mapping subproblem retains all residual terms
    involving the updated vertices, including those from adjacent simplices and
    any constrained boundary faces. The auxiliary and dual variables can still
    be updated globally. Since the repair region is usually much smaller than
    the entire mesh, this localization can considerably reduce the computational
    cost of Stage~II. Once all foldings have been removed, the resulting strictly
    feasible map is used to initialize Stage~III.
\end{remark}

\subsection{Stage III: Strictly feasible barrier minimization}
\label{subsec:stage-III-energy-minimization}

With folded and degenerate simplices removed in Stage~II, we now
improve the parameterization by minimizing $E_{\beta,\eta}$ in
\cref{eq:discrete_joint_energy} over $\calX_{\mathtt B}$, starting from
$f^{\mathrm{II}}$.
The barrier penalizes compression toward degeneracy as angular and
volumetric distortions are reduced. By \Cref{sec:DiscreteEnergy}, the
Euclidean gradient of $E_{\beta,\eta}$ is
\begin{equation*}
    \nabla_{\bff}E_{\beta,\eta}(f)
    =\left[
    \frac{2\beta}{n}\bigl(L_D-L_A(f)\bigr)
    +2(1-\beta)\bigl(L_{V_s}(f)-L_{V_r}(f)\bigr)
    -\eta L_B(f)
    \right]\bff.
\end{equation*}
Its projection onto the tangent space gives the Riemannian gradient.
% $\operatorname{grad}E_{\beta,\eta}(f)=\mathcal P_{\bff}\!\left(\nabla_{\bff}E_{\beta,\eta}(f)\right)$.

We use a preconditioned Riemannian L-BFGS method~\cite{AMS08,HGA15,NW06} with backtracking line search.
If the computed direction is not a descent direction, we use the
preconditioned negative gradient instead. A step is accepted when all
simplex Jacobians remain positive and both $E_\beta$ and $E_{\beta,\eta}$
satisfy the sufficient-decrease conditions. For ball parameterizations,
we also require every boundary $(n-1)$-simplex to retain its outward
orientation relative to $\bbS^{n-1}$.

The positivity test excludes invalid trial steps, while the barrier
influences the search direction through its gradient, which gives
increasing weight to simplices with small Jacobians. In our preliminary
experiments, minimizing $E_\beta$ without the barrier but retaining the
positivity test often led to repeated step reductions and eventual
line-search failure. The barrier helps avoid this difficulty by
discouraging further compression of nearly degenerate simplices,
making it easier to find feasible steps that reduce the energy.

For spherical parameterizations, we may also limit the vertex
displacements during line search so that successive iterates satisfy
the conditions of \Cref{lem:degree-preservation}. Their radial
normalizations are then homotopic. Hence, if the map entering Stage~III
has degree one after radial normalization, this degree is preserved
throughout the minimization.

\section{Numerical experiments}
\label{sec:numerical-experiments}

Numerical experiments are conducted on four classes of parameterization
problems: disk parameterizations of open triangular surfaces, spherical
parameterizations of closed triangular surfaces, ball parameterizations of
tetrahedral volumes, and $\bbS^3$ parameterizations of closed
three-dimensional hypersurfaces embedded in $\re^4$. The test set spans
intrinsic dimensions two and three and includes both ball and sphere targets.
{The numerical comparisons report global energies, local distortion
statistics, and bijectivity measures.}
For each $\tau\in\Sigma_n(\calM)$, let
$\sigma_1,\ldots,\sigma_n$ be singular values of
$\nabla_{\calM}f$, as defined in
\Cref{sec:ContinuousEnergy}, applied to the affine restriction $f|_\tau$.
{To quantify angular and volumetric distortions,} we define
simplexwise conformality distortion, conformal error, and volume-ratio error by $K_\tau(f)={1\over n}{\sum_{i=1}^n\sigma_i^2}/
{\left(\prod_{i=1}^n\sigma_i\right)^{2/n}}$, $d_{C,\tau}=K_\tau(f)-1$, $d_{V,\tau}=
\abs{{J_\tau(f)}/{\rmR_f}-1}$. Then the weighted mean and standard deviation are,
for $\star\in\{C,V\}$, $\mu_\star=\sum\limits_{\tau\in\Sigma_n(\calM)}d_{\star,\tau}{\abs{\tau}\over\mu(\calM)}$, $s_\star=\left[\sum\limits_{\tau\in\Sigma_n(\calM)}(d_{\star,\tau}-\mu_\star)^2{\abs{\tau}\over\mu(\calM)}\right]^{1/2}$.
% \[\mu_\star=
% \sum_{\tau\in\Sigma_n(\calM)}d_{\star,\tau}\cdot{\abs{\tau}}/{\mu(\calM)},
% \qquad
% s_\star=
% \left[
% \sum_{\tau\in\Sigma_n(\calM)}
% (d_{\star,\tau}-\mu_\star)^2\cdot{\abs{\tau}}/{\mu(\calM)}
% \right]^{1/2}.\]
We use ``fold'' to denote the number of simplices satisfying
$J_\tau(f)\leq 0$ and set
$J_{\min}=\min\limits_{\tau\in\Sigma_n(\calM)}J_\tau(f)$.
Consequently, \(\mathrm{fold}=0\), or equivalently \(J_{\min}>0\),
verifies {local injectivity on every simplex and excludes flipped or degenerate simplices}.
By \Cref{thm:global-bijectivity-certificate}, a spherical
parameterization is globally bijective if its radial normalization
is well defined, locally injective, and of degree one.
We verify local injectivity using the closed-vertex-star embedding
condition and check that the degree is one. For ball parameterizations,
global bijectivity follows from boundary bijectivity and $J_{\min}>0$.
All final parameterizations produced by the proposed method satisfy
these conditions.

\subsection{Parameterization of open surfaces onto the disk}

\begin{figure}
\centering
\includegraphics[width=0.90\textwidth]{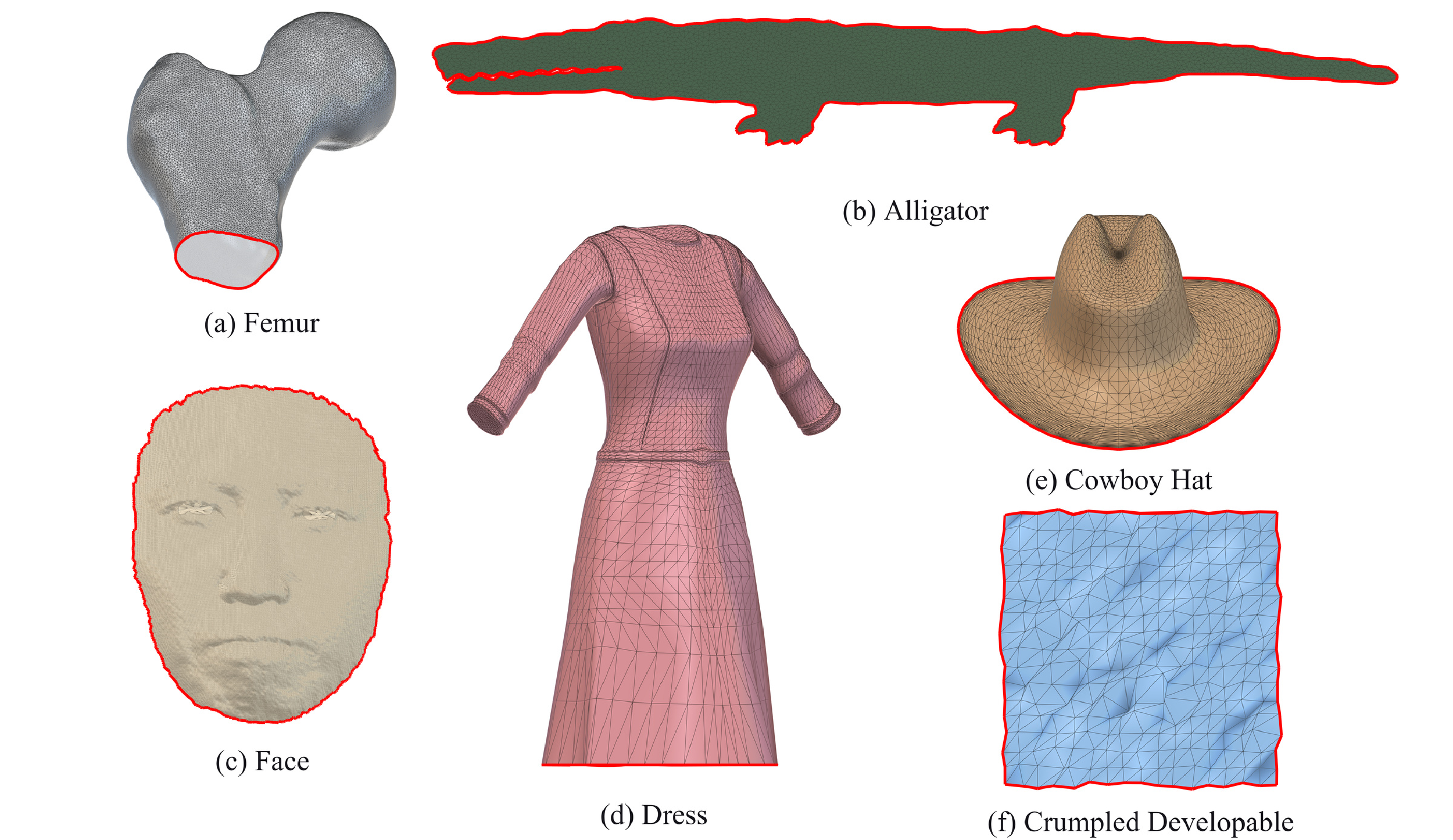}
\caption{Triangular benchmark models for parameterizations onto the disk. The red lines indicate the mesh boundaries.}
\label{fig:2d-models}
\end{figure}

The two-dimensional experiments use six publicly available triangular meshes
from the AIM@SHAPE Digital Shape Workbench~\cite{AIM},
{Keenan's 3D Model Repository~\cite{Keenan},
Sketchfab~\cite{sketchfab}, and the Stanford 3D
Scanning Repository~\cite{Stanford}}
{as shown in \cref{fig:2d-models}.} We set
$\beta=0.5$ and compare the proposed method with the energy-based
distortion-balancing parameterization method of Liu and
Yueh~\cite{LY25}\footnote{The implementation was downloaded from
\url{https://github.com/liu227857/BEM}.}, which seeks to achieve bijectivity
through a mean-value Laplacian whose positive weights, together with a
prescribed convex boundary, yield a convex-combination property.
Table~\ref{tab:disk2d} compares the two methods on the six open
triangular surfaces. Both methods produce fold-free maps for Femur and Cowboy Hat with comparable computation times. The method of Liu and Yueh does not produce usable results for Alligator or Crumpled Developable and retains 58 and 1 local folds for Face and Dress, respectively. In contrast, the proposed method maintains strictly positive signed local volume ratios for all six inputs.

\begin{table}
\centering
\caption{Comparison for disk parameterization of six open triangular surfaces.}
\label{tab:disk2d}
\setlength{\tabcolsep}{2.4pt}
\resizebox{0.85\textwidth}{!}{
\begin{tabular}{ccccccccccc}
\toprule
Model & Algorithm & $E_C$ & $E_V$
& $\mu_C$ & $s_C$ & $\mu_V$ & $s_V$
& $J_{\min}$ & Folds & Time (s) \\
\midrule

\multirow{2}{*}{Femur}
& Proposed ($\beta=0.5$)
& $1.44\times10^{0}$ & $1.06\times10^{1}$
& 0.107 & 0.280 & 1.478 & 1.097
& $6.64\times10^{-5}$ & 0 & 1.34 \\
& Liu \& Yueh~\cite{LY25}
& $2.87\times10^{0}$ & $2.87\times10^{0}$
& 0.389 & 0.634 & 0.810 & 0.509
& $5.38\times10^{-3}$ & 0 & 1.83 \\
\midrule

\multirow{2}{*}{Alligator}
& Proposed ($\beta=0.5$)
& $6.40\times10^{0}$ & $9.37\times10^{-1}$
& 3.117 & 11.478 & 0.444 & 0.318
& $1.78\times10^{-5}$ & 0 & 0.30 \\
& Liu \& Yueh~\cite{LY25}
& -- & -- & -- & -- & -- & --
& -- & -- & Failed \\
\midrule

\multirow{2}{*}{Face}
& Proposed ($\beta=0.5$)
& $1.17\times10^{-1}$ & $1.86\times10^{-1}$
& 0.079 & 1.863 & 0.178 & 0.165
& $1.37\times10^{-3}$ & 0 & 1.02 \\
& Liu \& Yueh~\cite{LY25}
& $7.81\times10^{-2}$ & $7.82\times10^{-2}$
& 0.250 & 63.316 & 0.108 & 0.115
& $-2.01\times10^{-1}$ & 58 & 4.74 \\
\midrule

\multirow{2}{*}{Dress}
& Proposed ($\beta=0.5$)
& $7.60\times10^{-2}$ & $1.02\times10^{1}$
& 0.021 & 0.077 & 1.326 & 1.224
& $3.77\times10^{-12}$ & 0 & 0.08 \\
& Liu \& Yueh~\cite{LY25}
& $2.23\times10^{0}$ & $2.23\times10^{0}$
& 0.429 & 0.345 & 0.777 & 0.333
& $-4.40\times10^{-3}$ & 1 & 0.20 \\
\midrule

\multirow{2}{*}{Cowboy Hat}
& Proposed ($\beta=0.5$)
& $1.89\times10^{-1}$ & $3.17\times10^{-1}$
& 0.060 & 0.051 & 0.257 & 0.189
& $1.78\times10^{-1}$ & 0 & 0.11 \\
& Liu \& Yueh~\cite{LY25}
& $2.01\times10^{-1}$ & $2.00\times10^{-1}$
& 0.059 & 0.055 & 0.224 & 0.118
& $4.82\times10^{-1}$ & 0 & 0.12 \\
\midrule

\multirow{2}{*}{Crumpled Dev.}
& Proposed ($\beta=0.5$)
& $1.62\times10^{-1}$ & $1.54\times10^{-1}$
& 1.577 & 24.397 & 0.146 & 0.166
& $7.74\times10^{-6}$ & 0 & 0.01 \\
& Liu \& Yueh~\cite{LY25}
& -- & -- & -- & -- & -- & --
& -- & -- & Failed \\
\bottomrule
\end{tabular}
}
\end{table}

\subsection{Surface and volumetric parameterizations in $\re^3$}

\begin{figure}
\centering
\includegraphics[width=0.95\textwidth]{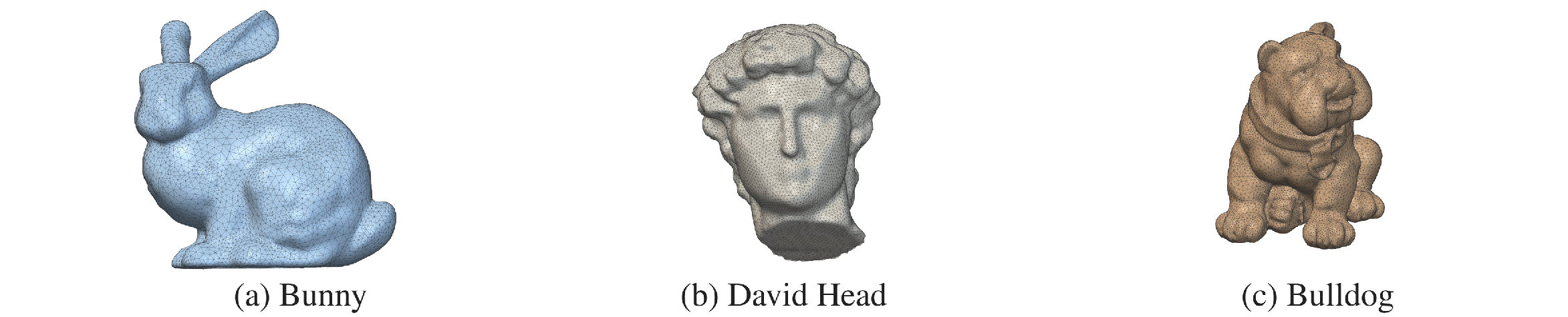}
\caption{Tetrahedral benchmark models for parameterizations in $\re^3$.}
\label{fig:3d-models}
\end{figure}

For the volumetric experiments, we use the publicly available tetrahedral meshes of Bunny~\cite{Stanford}, David Head~\cite{stanford-david}, and Bulldog~\cite{vcg-models}, shown in \cref{fig:3d-models}. Their boundary triangular meshes are extracted for the surface experiments. The comparison methods are spherical conformal energy minimization (CEM), spherical stretch energy minimization (SEM), and volumetric stretch energy minimization (VSEM), as proposed in~\cite{MHTL19}. For the surface experiments, we compare the proposed method with $\beta=0$ against SEM and the proposed method with $\beta=1$ against CEM. For the volumetric experiments, we compare the proposed method with $\beta=0$ against VSEM. The numerical results are described below.

\subsubsection{Boundary-surface parameterization onto $\bbS^2$}

\begin{figure}
\centering
\includegraphics[width=0.9\textwidth]{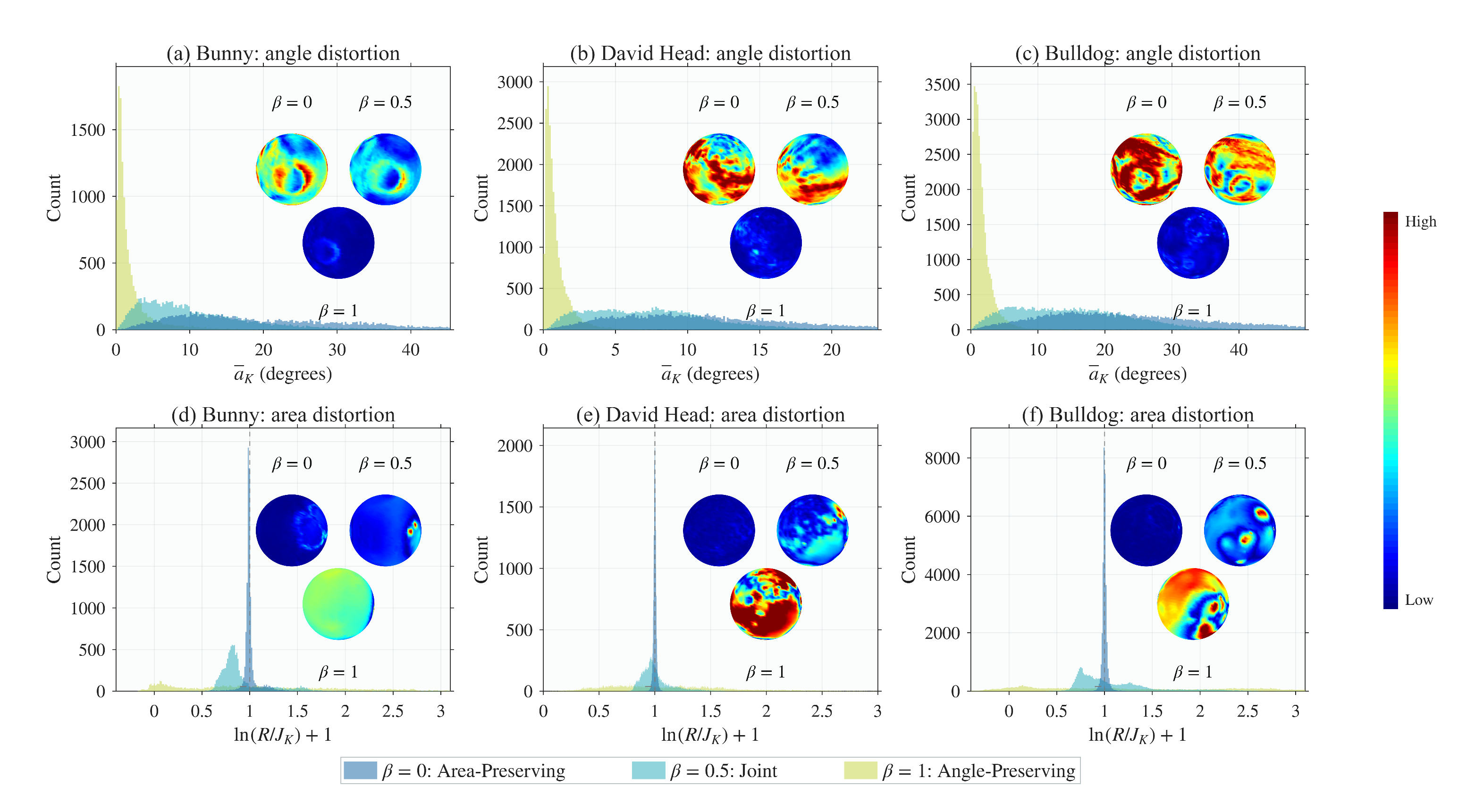}
\caption{Distributions of local surface distortion for Bunny, DavidHead, and Bulldog. The color bar characterizes the level of local distortion.}

\label{fig:surface-hist}
\end{figure}

Each histogram sample in \cref{fig:surface-hist} corresponds to one triangle
$K$. Let $\theta_{K,j}$ and $\widehat\theta_{K,j}$ be the $j$th interior angle
before and after parameterization, respectively, and let
$a_{K,j}=\widehat\theta_{K,j}-\theta_{K,j}$. The top row plots the mean
absolute interior-angle error
$\bar a_K=\frac13\sum_{j=1}^3\abs{a_{K,j}}$. The bottom row plots
$q_K=\ln(\rmR_f/J_K(f))+1$, for which $q_K=1$ means that the local area
ratio equals the global ratio. 
For a fixed model, the three values
of $\beta$ use a common color scale to measure the angular or volumetric distortion.

As shown in \Cref{tab:surface3d}, all nine spherical parameterizations
computed by the proposed method are fold-free. In contrast, SEM retains
four folds for Bunny, and CEM retains one fold for Bulldog. The energy
values exhibit a trade-off between $E_C$ and $E_V$, with $\beta=0.5$
balancing the two objectives. \Cref{fig:surface-hist} shows the
corresponding local angle and area distortions.

\begin{table}
\centering
\caption{Comparison for parameterization of closed two-dimensional surfaces onto $\bbS^2$.}
\label{tab:surface3d}
\scriptsize
\setlength{\tabcolsep}{2.4pt}
\resizebox{0.85\textwidth}{!}{
\begin{tabular}{ccccccccccc}
\toprule
Model & Algorithm & $E_C$ & $E_V$
& $\mu_C$ & $s_C$ & $\mu_V$ & $s_V$
& $J_{\min}$ & Folds & Time (s) \\
\midrule

\multirow{5}{*}{Bunny}
& Proposed ($\beta=0$)
& $9.28\times10^{0}$ & $1.56\times10^{-1}$
& 0.783 & 1.914 & 0.048 & 0.100
& $3.94\times10^{-2}$ & 0 & 1.65 \\
& SEM~\cite{MHTL19}
& $9.23\times10^{0}$ & $1.59\times10^{-1}$
& 0.778 & 1.815 & 0.049 & 0.101
& $-2.65\times10^{-1}$ & 4 & 0.14 \\
\cmidrule(r){2-11}
& Proposed ($\beta=0.5$)
& $1.26\times10^{0}$ & $2.08\times10^{0}$
& 0.093 & 0.131 & 0.314 & 0.260
& $1.97\times10^{-5}$ & 0 & 0.98 \\
\cmidrule(r){2-11}
& Proposed ($\beta=1$)
& $9.04\times10^{-3}$ & $1.05\times10^{1}$
& 0.006 & 0.018 & 0.774 & 0.487
& $9.32\times10^{-8}$ & 0 & 0.63 \\
& CEM~\cite{MHTL19}
& $1.35\times10^{-2}$ & $1.05\times10^{1}$
& 0.006 & 0.013 & 0.663 & 0.629
& $6.66\times10^{-8}$ & 0 & 0.09 \\
\midrule

\multirow{5}{*}{David Head}
& Proposed ($\beta=0$)
& $1.22\times10^{0}$ & $2.28\times10^{-3}$
& 0.097 & 0.115 & 0.009 & 0.010
& $5.95\times10^{-1}$ & 0 & 0.95 \\
& SEM~\cite{MHTL19}
& $1.21\times10^{0}$ & $1.36\times10^{-3}$
& 0.096 & 0.116 & 0.007 & 0.008
& $8.63\times10^{-1}$ & 0 & 1.48 \\
\cmidrule(r){2-11}
& Proposed ($\beta=0.5$)
& $4.73\times10^{-1}$ & $2.01\times10^{-1}$
& 0.036 & 0.038 & 0.095 & 0.084
& $4.50\times10^{-1}$ & 0 & 1.83 \\
\cmidrule(r){2-11}
& Proposed ($\beta=1$)
& $7.87\times10^{-3}$ & $2.70\times10^{0}$
& 0.001 & 0.002 & 0.388 & 0.254
& $7.89\times10^{-2}$ & 0 & 0.30 \\
& CEM~\cite{MHTL19}
& $7.87\times10^{-3}$ & $2.69\times10^{0}$
& 0.001 & 0.002 & 0.387 & 0.255
& $8.04\times10^{-2}$ & 0 & 0.15 \\
\midrule

\multirow{5}{*}{Bulldog}
& Proposed ($\beta=0$)
& $5.53\times10^{0}$ & $6.14\times10^{-3}$
& 0.441 & 0.457 & 0.014 & 0.018
& $9.16\times10^{-2}$ & 0 & 2.98 \\
& SEM~\cite{MHTL19}
& $5.57\times10^{0}$ & $5.70\times10^{-3}$
& 0.444 & 0.466 & 0.014 & 0.016
& $1.41\times10^{-1}$ & 0 & 1.65 \\
\cmidrule(r){2-11}
& Proposed ($\beta=0.5$)
& $2.28\times10^{0}$ & $1.08\times10^{0}$
& 0.160 & 0.171 & 0.246 & 0.159
& $2.22\times10^{-3}$ & 0 & 4.94 \\
\cmidrule(r){2-11}
& Proposed ($\beta=1$)
& $1.45\times10^{-2}$ & $1.16\times10^{1}$
& 0.005 & 2.162 & 0.813 & 0.515
& $1.06\times10^{-6}$ & 0 & 0.79 \\
& CEM~\cite{MHTL19}
& $1.54\times10^{-2}$ & $1.33\times10^{1}$
& 0.002 & 0.140 & 0.831 & 0.611
& $-4.18\times10^{-3}$ & 1 & 0.24 \\
\bottomrule
\end{tabular}
}
\end{table}

\subsubsection{Volumetric parameterization onto $\bbB^3$}

\begin{figure}
\centering
\includegraphics[width=0.9\textwidth]{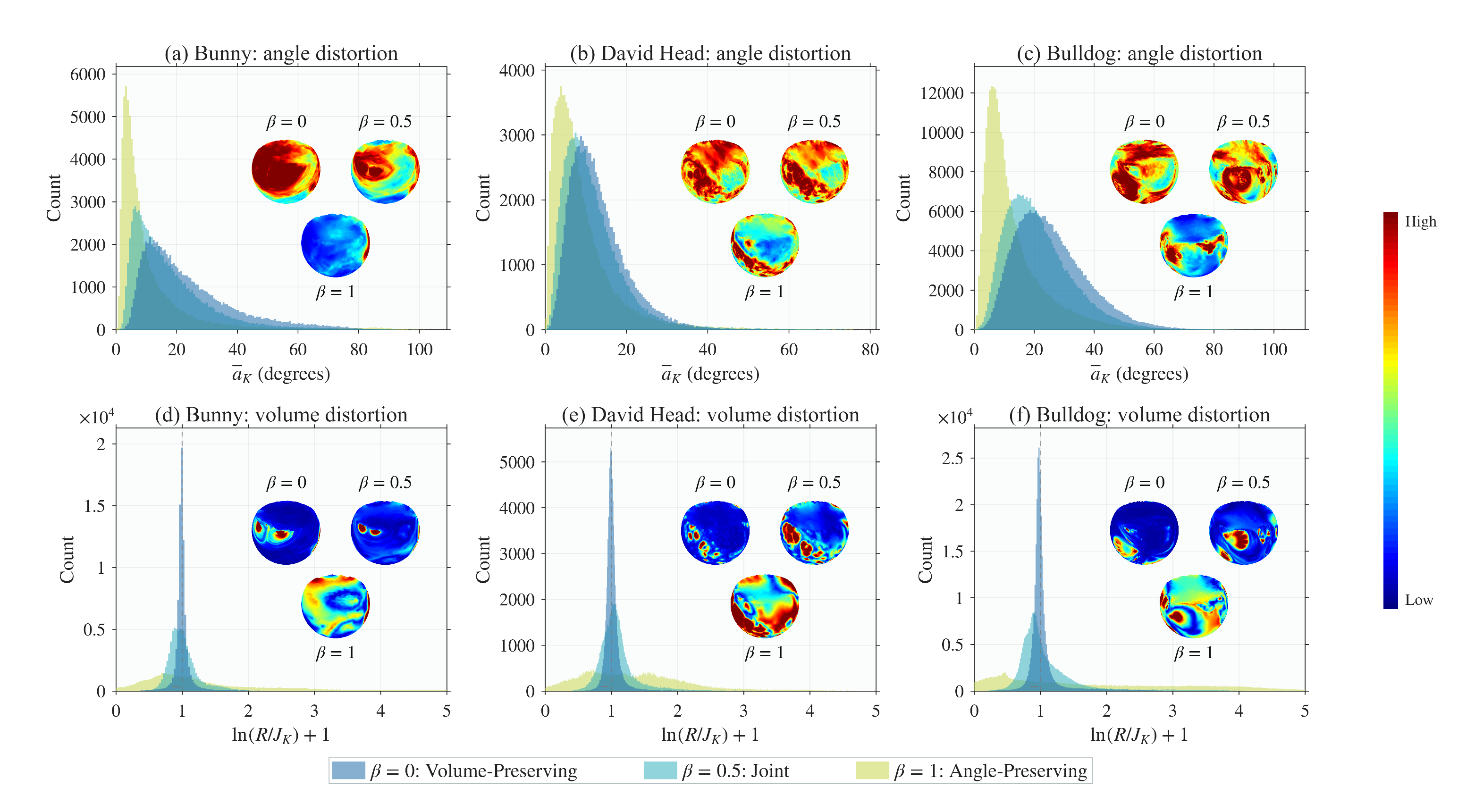}
\caption{Distributions of local volumetric distortion for Bunny, DavidHead,
and Bulldog. The color bar characterizes the level of local distortion.}
\label{fig:volume-hist}
\end{figure}

\Cref{fig:volume-hist} {follows the same statistical conventions as}
the surface figure, but each histogram sample corresponds to one tetrahedron. Let
$a_{K,j}$ be the change in the $j$th dihedral angle under parameterization.
The top row plots the mean absolute error over the six dihedral angles,
$\bar a_K=\frac16\sum_{j=1}^6\abs{a_{K,j}}$. The bottom row again plots
$q_K=\ln(\rmR_f/J_K(f))+1$, where $q_K=1$ means that the local volume ratio
equals the global ratio.
{For a fixed model, the three values of $\beta$ use a common color
scale to measure the angular or volumetric distortion.}

The proposed method also produces fold-free volumetric parameterizations
for all three models and all three choices of $\beta$
(\Cref{tab:volume3d}). VSEM produces 106, 49, and 10 folds for Bulldog,
Bunny, and DavidHead, respectively. The energy values again show the
trade-off between $E_C$ and $E_V$, with $\beta=0.5$ providing an
intermediate balance.

\begin{table}
\centering
\caption{Comparison for volumetric parameterization onto
$\bbB^3$.}
\label{tab:volume3d}
\scriptsize
\setlength{\tabcolsep}{2.4pt}
\resizebox{0.85\textwidth}{!}{
\begin{tabular}{ccccccccccc}
\toprule
Model & Algorithm & $E_C$ & $E_V$
& $\mu_C$ & $s_C$ & $\mu_V$ & $s_V$
& $J_{\min}$ & Folds & Time (s) \\
\midrule

\multirow{4}{*}{Bulldog}
& Proposed ($\beta=0$)
& $2.50\times10^{0}$ & $8.44\times10^{-2}$
& 1.137 & 137.643 & 0.075 & 0.121
& $9.48\times10^{-7}$ & 0 & 55.29 \\
& VSEM~\cite{MHTL19}
& $2.94\times10^{0}$ & $2.30\times10^{-2}$
& 0.714 & 1.265 & 0.046 & 0.059
& $-1.41\times10^{0}$ & 106 & 32.74 \\
\cmidrule(r){2-11}
& Proposed ($\beta=0.5$)
& $1.62\times10^{0}$ & $3.21\times10^{-1}$
& 0.420 & 82.692 & 0.208 & 0.182
& $1.59\times10^{-7}$ & 0 & 41.99 \\
\cmidrule(r){2-11}
& Proposed ($\beta=1$)
& $3.37\times10^{-1}$ & $3.97\times10^{0}$
& 0.152 & 2.403 & 0.809 & 0.544
& $7.81\times10^{-9}$ & 0 & 64.34 \\
\midrule

\multirow{4}{*}{Bunny}
& Proposed ($\beta=0$)
& $2.66\times10^{0}$ & $1.26\times10^{-1}$
& 2.049 & 206.565 & 0.066 & 0.161
& $1.04\times10^{-5}$ & 0 & 17.30 \\
& VSEM~\cite{MHTL19}
& $3.07\times10^{0}$ & $5.16\times10^{-2}$
& 0.802 & 3.996 & 0.055 & 0.097
& $-3.54\times10^{0}$ & 49 & 9.38 \\
\cmidrule(r){2-11}
& Proposed ($\beta=0.5$)
& $1.11\times10^{0}$ & $2.36\times10^{-1}$
& 0.638 & 6.425 & 0.160 & 0.176
& $4.46\times10^{-11}$ & 0 & 8.17 \\
\cmidrule(r){2-11}
& {Proposed ($\beta=1$)}
& $2.34\times10^{-1}$ & $2.23\times10^{0}$
& 3.959 & 70.106 & 0.557 & 0.473
& $3.07\times10^{-17}$ & 0 & 4.80 \\
\midrule

\multirow{4}{*}{David Head}
& Proposed ($\beta=0$)
& $8.68\times10^{-1}$ & $2.62\times10^{-2}$
& 0.213 & 1.106 & 0.049 & 0.062
& $9.04\times10^{-5}$ & 0 & 6.50 \\
& VSEM~\cite{MHTL19}
& $9.24\times10^{-1}$ & $1.83\times10^{-2}$
& 0.223 & 0.294 & 0.044 & 0.050
& $-1.42\times10^{-1}$ & 10 & 7.49 \\
\cmidrule(r){2-11}
& Proposed ($\beta=0.5$)
& $5.82\times10^{-1}$ & $1.09\times10^{-1}$
& 0.148 & 2.282 & 0.108 & 0.119
& $1.02\times10^{-5}$ & 0 & 4.79 \\
\cmidrule(r){2-11}
& Proposed ($\beta=1$)
& $2.44\times10^{-1}$ & $1.11\times10^{0}$
& 0.075 & 0.234 & 0.415 & 0.307
& $1.82\times10^{-4}$ & 0 & 2.74 \\
\bottomrule
\end{tabular}
}
\end{table}

\subsection{Hypersurface parameterization onto $\bbS^3$}

We test the proposed method on two closed three-dimensional hypersurfaces embedded in \(\mathbb{R}^4\). Both models consist of 378 vertices and 2114 3-simplices. We consider \(\beta\in\{0,0.5,1\}\) and use the \(n\)-VSEM method in~\cite{ZHTL25} for comparison. The two test manifolds $\calM_i\subseteq\re^4$ ($i=1,2$) are defined as
\[\calM_1=\{Ds:\,s\in\mathbb{S}^3,\, D=\diag(0.7,\,0.9,\,1.1,\,1.3)\},\ \calM_2=\{\rho(s)s:\,s\in\mathbb{S}^3,\,\rho(s)=0.5+s_4^2\}.\]

\begin{table}
\centering
\caption{Comparison for parameterization of three-dimensional
hypersurfaces in $\re^4$ onto $\bbS^3$.}
\label{tab:surface4d}
\scriptsize
\setlength{\tabcolsep}{2.4pt}
\renewcommand{\arraystretch}{1.08}
\resizebox{0.85\textwidth}{!}{
\begin{tabular}{ccccccccccc}
\toprule
Model & Algorithm & $E_C$ & $E_V$
& $\mu_C$ & $s_C$ & $\mu_V$ & $s_V$
& $J_{\min}$ & Folds & Time (s) \\
\midrule

\multirow{4}{*}{Ellipsoid}
& Proposed ($\beta=0$)
& $1.11\times10^{0}$ & $7.64\times10^{-3}$
& 0.060 & 0.043 & 0.017 & 0.014
& $8.14\times10^{-1}$ & 0 & 11.00 \\
& $n$-VSEM~\cite{ZHTL25}
& $3.59\times10^{0}$ & $3.16\times10^{0}$
& 0.198 & 0.159 & 0.367 & 0.251
& $1.01\times10^{-1}$ & 0 & 0.09 \\
\cmidrule(r){2-11}
& Proposed ($\beta=0.5$)
& $4.54\times10^{-1}$ & $1.35\times10^{-1}$
& 0.024 & 0.009 & 0.076 & 0.050
& $7.21\times10^{-1}$ & 0 & 3.53 \\
\cmidrule(r){2-11}
& Proposed ($\beta=1$)
& $1.16\times10^{-1}$ & $1.67\times10^{0}$
& 0.007 & 0.005 & 0.264 & 0.182
& $3.91\times10^{-1}$ & 0 & 3.08 \\
\midrule

\multirow{4}{*}{Peanut}
& Proposed ($\beta=0$)
& $6.27\times10^{0}$ & $5.14\times10^{-1}$
& 0.356 & 0.564 & 0.143 & 0.132
& $2.24\times10^{-1}$ & 0 & 13.70 \\
& $n$-VSEM~\cite{ZHTL25}
& $1.68\times10^{1}$ & $1.26\times10^{1}$
& 1.589 & 1.540 & 0.809 & 1.054
& $-3.77\times10^{0}$ & 141 & 0.03 \\
\cmidrule(r){2-11}
& Proposed ($\beta=0.5$)
& $3.49\times10^{0}$ & $6.17\times10^{-1}$
& 0.176 & 0.192 & 0.173 & 0.114
& $5.77\times10^{-1}$ & 0 & 6.31 \\
\cmidrule(r){2-11}
& Proposed ($\beta=1$)
& $1.69\times10^{-1}$ & $8.96\times10^{0}$
& 0.064 & 0.081 & 0.934 & 0.296
& $1.92\times10^{-3}$ & 0 & 3.72 \\
\bottomrule
\end{tabular}
}
\end{table}

All six parameterizations computed by the proposed method are fold-free,
as shown in \Cref{tab:surface4d}. Both models show the same energy
trade-off: $\beta=0$ gives lower $E_V$, $\beta=1$ gives lower $E_C$,
and $\beta=0.5$ lies between the two endpoints.
Although the proposed method requires more CPU time than $n$-VSEM, it produces fold-free results and attains lower energy than $n$-VSEM on both models in the volume-preserving case ($\beta=0$).

\section{Application to medical imaging}

This section applies the proposed parameterization method to brain-tumor label
transfer between a patient-specific brain domain and a canonical cube.
Patient-specific brain domains have irregular boundaries and vary substantially
in shape and spatial extent, whereas convolutional neural networks and other
medical-image analysis methods commonly require fixed-size regular
three-dimensional arrays.
{Resampling tumor labels on a canonical cube provides a
fixed-size representation. Foldings may introduce ambiguous sampling
correspondences and affect label transfer.}
The experiment is designed to determine whether enforcing bijectivity improves
the geometric validity of the parameterization and
{the accuracy of label transfer and native-grid reconstruction.}
{The brain and cube parameterizations are represented in the
same unit-ball coordinate system.}
The complete bijective method is compared with a non-bijective ablation under
the same data, discretization, and resampling settings. Both parameterizations
are computed independently of the tumor labels, which are used only to evaluate
the resulting fixed maps.

Let $\calM\subset\re^3$ be the brain domain represented by a patient-specific
tetrahedral mesh,
{and let $\mathcal C=[-1,1]^3$ be the canonical cube.
Denote the piecewise-affine maps by $f:\calM\to\Omega_f$ and
$g:\mathcal C\to\Omega_g$, with boundary vertex images on the unit sphere.
The two polyhedral images differ slightly near their boundaries.
The gliomas considered here lie in the brain interior, so these differences
do not affect the present experiments. For points mapped into
$\Omega_f\cap\Omega_g$, the composition $g^{-1}\circ f$ gives the
corresponding cube coordinates.}

\begin{figure}
\centering
\includegraphics[width=0.9\textwidth]
{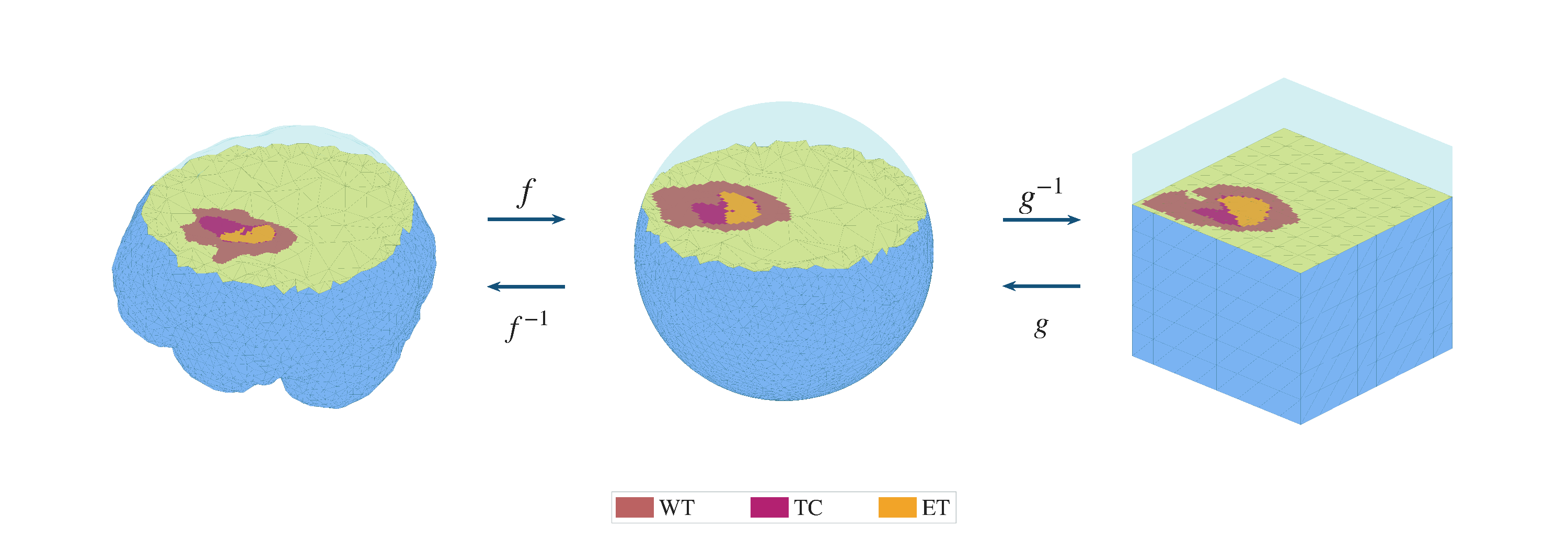}
\caption{{Brain-tumor label transfer between the
patient-specific brain domain and the canonical cube through their
ball parameterizations.}}
\label{fig:brain-parameterization}
\end{figure}

The parameterization and label-transfer pipeline involves several
discrete representations, as illustrated in \cref{fig:brain-parameterization}.
A patient-specific tetrahedral mesh is used to compute
$f$, while the original NIfTI voxel grid stores the brain image and tumor
labels. Each coordinate direction of $\mathcal C$ is divided into $12$ equal
intervals, and every small cube is split consistently into six tetrahedra. The
resulting $10{,}368$ tetrahedra are used to compute $g$. Independently, a
regular $128^3$ voxel grid on $\mathcal C$ stores the resampled labels and does
not participate in the parameterization.

Let $A\subset\calM$ be a tumor region on the original NIfTI grid, let $C_j$ be
the $j$th voxel of the canonical cube, and let $w_j$ be its center.
{To resample the labels, we locate $g(w_j)$ in the mapped
brain mesh and use the barycentric coordinates in the containing
tetrahedron to obtain the corresponding native point $\bfx_j$.
The original NIfTI label is then sampled at $\bfx_j$ by nearest-neighbor
interpolation. If $\bfx_j\in A$, the cube voxel $C_j$ is marked as tumor.
The union of all marked voxels, denoted by $B\subset\mathcal C$, gives the
resampled tumor region on the regular $128^3$ grid.}

For the return transfer, the algorithm uses the points $\bfx_j$ stored
during the forward resampling and scatters
{their associated labels, including background labels,}
back to the original brain grid.
If several cube samples enter the same brain voxel, the sample closest to
the voxel center is retained. A voxel receiving no direct sample is filled from
the nearest mapped sample in physical space.
The reconstructed region, $\widehat A$, lies on the same NIfTI grid and
in the same physical coordinate system as $A$.

The reconstruction accuracy is evaluated by comparing $\widehat A$ directly
with $A$ on the original NIfTI grid. Let $\partial A$ and
$\partial\widehat A$ denote their sets of boundary-voxel centers, and define
$d(\bfx,S)=\min_{\bfy\in S}\norm{\bfx-\bfy}$. Region overlap is measured by the Dice
similarity coefficient (DSC), while boundary agreement is measured by the
average symmetric surface distance (ASSD):
\begin{equation*}
  \operatorname{DSC}(A,\widehat A)=
  \frac{2\abs{A\cap\widehat A}}
       {\abs{A}+\abs{\widehat A}},\qquad
  \operatorname{ASSD}(A,\widehat A)=
  \frac{\sum_{\bfx\in\partial A}d(\bfx,\partial\widehat A)
    +\sum_{\bfy\in\partial\widehat A}d(\bfy,\partial A)}
  {\abs{\partial A}+\abs{\partial\widehat A}}.
\end{equation*}
We also record the numbers of false-positive and false-negative voxels, i.e.,
$\mathrm{FP}=\#(\widehat A\setminus A)$,
$\mathrm{FN}=\#(A\setminus\widehat A)$.
A larger DSC and smaller ASSD, FP, and FN indicate better reconstruction
accuracy.

The experiment contains six glioma cases from the BraTS 2020 training dataset
\cite{Bakas2017TCGA,Bakas2018BRATS,Menze2015BRATS}, denoted by C1--C6.
Whole tumor (WT), tumor core (TC), and enhancing tumor (ET) are evaluated
in each case. The patient-specific brain meshes contain
$5{,}162$--$5{,}924$ vertices and $20{,}760$--$24{,}145$ tetrahedra.
We use $\beta=0.5$ for both the brain and cube parameterizations.

The experiment follows a strictly paired design. The bijective and
non-bijective arms use the same brain mesh, original image, regular $128^3$
cube grid, and bijective cube map $g$; they differ only in the construction of
the brain map $f$. The bijective arm uses the complete proposed three-stage
algorithm, with fold repair performed by ADMM--PNCG in Stage~II and the
joint energy with a Jacobian barrier, $E_{\beta,\eta}$, minimized in the
strictly feasible region in Stage~III.
The a posteriori conditions in
\Cref{thm:global-bijectivity-certificate} were also verified for the
brain-to-ball and cube-to-ball maps used in the bijective arm.
The non-bijective ablation omits all surface and volumetric fold repair and
minimizes the same joint energy $E_{\beta}$ without the Jacobian barrier.

For each boundary triangle, we check whether its orientation is
preserved, while for each tetrahedron $\tau$, we compute the signed local
volume ratio $J_\tau(f)$ defined in
\cref{eq:signed_local_volume_ratio}. A boundary triangle with reversed or
degenerate orientation and a tetrahedron with $J_\tau(f)\leq 0$ are counted
as a surface fold and a volumetric fold, respectively.
For folded tetrahedra, we also accumulate their physical volumes
in the source brain domain. The case-wise geometric results are given in
Table~\ref{tab:medical-casewise}. The bijective maps contain no surface or
volumetric folds, and their values of $J_{\min}$ range from
$3.40\times10^{-6}$ to $3.99\times10^{-5}$.
All six non-bijective maps contain folds, totaling $517$ folded
boundary faces and $1{,}031$ folded tetrahedra. Their mean surface- and
volume-fold rates are $1.186\%$ and $0.774\%$, respectively. The folded
tetrahedra have a mean source-domain physical volume of
$1{,}155\,\mathrm{mm}^3$, approximately $0.0694\%$ of the brain volume.

\begin{table}
  \centering
  \caption{{Case-wise geometric validity.}
  B and N denote the bijective and non-bijective arms, respectively;
  SF and VF are the numbers of folded boundary faces and tetrahedra;
  and $V^-$ is the source-domain physical volume (in $\mathrm{mm}^3$) of the
  tetrahedra satisfying $J_\tau(f)\leq0$. The corresponding rates are shown
  in parentheses.}
  \label{tab:medical-casewise}
  \scriptsize
  \setlength{\tabcolsep}{1.4pt}
  \begin{tabular}{@{}ccccccccc@{}}
    \toprule
    Case & SF (B) & SF (N) & VF (B) & VF (N)
      & $V^-$ (B) & $V^-$ (N) & $J_{\min}$ (B) & $J_{\min}$ (N)\\
    \midrule
    C1 & \textbf{0 (0.00\%)} & 123 (1.75\%) & \textbf{0 (0.00\%)} & 230 (1.09\%)
      & \textbf{0.0 (0.00\%)} & 969.1 (0.06\%) & $\boldsymbol{6.03\times10^{-6}}$ & $-1.731$\\
    \addlinespace[2pt]
    C2 & \textbf{0 (0.00\%)} & 65 (0.83\%) & \textbf{0 (0.00\%)} & 147 (0.61\%)
      & \textbf{0.0 (0.00\%)} & 944.2 (0.05\%) & $\boldsymbol{3.11\times10^{-5}}$ & $-1.253$\\
    \addlinespace[2pt]
    C3 & \textbf{0 (0.00\%)} & 137 (1.88\%) & \textbf{0 (0.00\%)} & 230 (1.03\%)
      & \textbf{0.0 (0.00\%)} & 1318.3 (0.08\%) & $\boldsymbol{3.99\times10^{-5}}$ & $-1.922$\\
    \addlinespace[2pt]
    C4 & \textbf{0 (0.00\%)} & 4 (0.05\%) & \textbf{0 (0.00\%)} & 75 (0.31\%)
      & \textbf{0.0 (0.00\%)} & 890.3 (0.05\%) & $\boldsymbol{1.73\times10^{-5}}$ & $-1.242$\\
    \addlinespace[2pt]
    C5 & \textbf{0 (0.00\%)} & 119 (1.71\%) & \textbf{0 (0.00\%)} & 211 (1.02\%)
      & \textbf{0.0 (0.00\%)} & 955.0 (0.06\%) & $\boldsymbol{2.63\times10^{-5}}$ & $-1.905$\\
    \addlinespace[2pt]
    C6 & \textbf{0 (0.00\%)} & 69 (0.89\%) & \textbf{0 (0.00\%)} & 138 (0.58\%)
      & \textbf{0.0 (0.00\%)} & 1853.0 (0.11\%) & $\boldsymbol{3.40\times10^{-6}}$ & $-2.240$\\
    \bottomrule
  \end{tabular}
\end{table}

The label-reconstruction results are reported in
Table~\ref{tab:medical-label-results}. The bijective arm achieves a higher DSC
and a lower ASSD in all $18$ paired comparisons. Its mean DSC values for WT,
TC, and ET are $99.54\%$, $99.43\%$, and $98.72\%$, compared with $99.39\%$,
$99.10\%$, and $98.20\%$ for the non-bijective arm. The corresponding gains
are $0.155$, $0.329$, and $0.516$ percentage points. The overall mean DSC
increases from $98.90\%$ to $99.23\%$, while the total number of false-positive
and false-negative voxels decreases from $17{,}548$ to $14{,}021$, a reduction
of $20.1\%$. The improvement is most pronounced for the smaller ET region,
which is more sensitive to local non-injectivity and discrete resampling.

\begin{table}
  \centering
  \caption{Case-wise label-reconstruction results for the three tumor regions.
  $\Delta\mathrm{DSC}=\mathrm{DSC}_{\mathrm B}-\mathrm{DSC}_{\mathrm N}$,
  and
  $\Delta\mathrm{ASSD}=\mathrm{ASSD}_{\mathrm N}-\mathrm{ASSD}_{\mathrm B}$.}
  \label{tab:medical-label-results}
  \scriptsize
  \setlength{\tabcolsep}{1.6pt}
  \begin{tabular}{@{}cccccccccc@{}}
    \toprule
    Case & Label & DSC (B) & DSC (N) & $\Delta$DSC
      & ASSD (B) & ASSD (N) & $\Delta$ASSD
      & FP/FN (B) & FP/FN (N)\\
    \midrule
    \multirow{3}{*}{C1}
      & WT & \textbf{99.71\%} & 99.57\% & 0.136\%
      & \textbf{0.0304} & 0.0420 & 0.0116
      & \textbf{221/211} & 308/324\\
      & TC & \textbf{99.45\%} & 99.00\% & 0.457\%
      & \textbf{0.0299} & 0.0506 & 0.0207
      & \textbf{56/53} & 100/100\\
      & ET & \textbf{99.31\%} & 98.84\% & 0.473\%
      & \textbf{0.0345} & 0.0518 & 0.0173
      & \textbf{39/37} & 58/70\\
    \addlinespace[2pt]
    \multirow{3}{*}{C2}
      & WT & \textbf{99.37\%} & 99.27\% & 0.100\%
      & \textbf{0.0617} & 0.0786 & 0.0170
      & \textbf{668/630} & 734/769\\
      & TC & \textbf{99.05\%} & 98.87\% & 0.182\%
      & \textbf{0.0710} & 0.0943 & 0.0234
      & \textbf{487/466} & 566/570\\
      & ET & \textbf{97.64\%} & 97.09\% & 0.549\%
      & \textbf{0.0621} & 0.0819 & 0.0198
      & \textbf{745/718} & 892/910\\
    \addlinespace[2pt]
    \multirow{3}{*}{C3}
      & WT & \textbf{99.76\%} & 99.52\% & 0.236\%
      & \textbf{0.0247} & 0.0518 & 0.0270
      & \textbf{160/136} & 277/308\\
      & TC & \textbf{99.77\%} & 99.00\% & 0.765\%
      & \textbf{0.0120} & 0.0499 & 0.0379
      & \textbf{14/10} & 53/50\\
      & ET & \textbf{99.82\%} & 98.96\% & 0.862\%
      & \textbf{0.0104} & 0.0481 & 0.0377
      & \textbf{8/9} & 41/56\\
    \addlinespace[2pt]
    \multirow{3}{*}{C4}
      & WT & \textbf{99.03\%} & 99.00\% & 0.029\%
      & \textbf{0.0663} & 0.0688 & 0.0025
      & \textbf{1934/1798} & 1940/1902\\
      & TC & \textbf{98.83\%} & 98.73\% & 0.098\%
      & \textbf{0.0700} & 0.0731 & 0.0031
      & \textbf{429/428} & 474/455\\
      & ET & \textbf{96.95\%} & 96.67\% & 0.276\%
      & \textbf{0.0583} & 0.0671 & 0.0088
      & \textbf{575/578} & 637/621\\
    \addlinespace[2pt]
    \multirow{3}{*}{C5}
      & WT & \textbf{99.86\%} & 99.57\% & 0.295\%
      & \textbf{0.0132} & 0.0417 & 0.0285
      & \textbf{104/92} & 320/295\\
      & TC & \textbf{99.87\%} & 99.46\% & 0.409\%
      & \textbf{0.0106} & 0.0389 & 0.0283
      & \textbf{36/36} & 161/142\\
      & ET & \textbf{99.79\%} & 99.10\% & 0.689\%
      & \textbf{0.0096} & 0.0368 & 0.0272
      & \textbf{48/48} & 214/191\\
    \addlinespace[2pt]
    \multirow{3}{*}{C6}
      & WT & \textbf{99.52\%} & 99.39\% & 0.133\%
      & \textbf{0.0455} & 0.0583 & 0.0129
      & \textbf{787/798} & 994/1029\\
      & TC & \textbf{99.62\%} & 99.56\% & 0.063\%
      & \textbf{0.0425} & 0.0494 & 0.0068
      & \textbf{277/297} & 334/335\\
      & ET & \textbf{98.82\%} & 98.57\% & 0.250\%
      & \textbf{0.0356} & 0.0424 & 0.0068
      & \textbf{541/547} & 653/665\\
    \midrule
    Mean & -- & \textbf{99.23\%} & 98.90\% & 0.333\%
      & \textbf{0.0382} & 0.0570 & 0.0187
      & \textbf{396.1/382.9} & 486.4/488.4\\
    \bottomrule
  \end{tabular}
\end{table}

{The complete method produces fold-free brain parameterizations
and yields a consistent improvement in tumor-label reconstruction
in every paired comparison under the same sampling and reconstruction
procedure.}
\section{Conclusion}
We developed a unified framework for bijective joint angle- and volume-preserving parameterizations of discrete $n$-manifolds.
The continuous formulation combines conformal and volume-preserving energies to provide an adjustable trade-off between local anisotropy and volumetric distortion.
At the discrete level, the constituent energies admit a unified cotangent Laplacian-type gradient representation, enabling sparse $n$-dimensional computation.
An ADMM--PNCG folding-repair stage followed by strictly feasible barrier
minimization {maintains local injectivity by preventing} element inversion and collapse.
{The topological analysis gives sufficient conditions for discrete spherical
maps to have degree one and establishes
global-bijectivity criteria for spherical and ball parameterizations.} Experiments on disk, sphere, ball, and $\bbS^3$ targets produced fold-free maps with effective distortion control and favorable comparisons with existing methods. {The a posteriori checks in \Cref{thm:global-bijectivity-certificate} further confirmed global bijectivity for all reported final parameterizations produced by the proposed method.}
In the brain-to-cube application, enforcing bijectivity eliminated all surface and volumetric folds and improved forward--inverse consistency and tumor-label reconstruction in every paired comparison.
These results demonstrate the potential of the framework as a general tool for geometry processing and medical-image analysis.

\section*{Acknowledgments}
We thank Tianhe-2 and the Big Data Computing Center in Southeast University, China, for their support of our use of their computing resources. The authors used AI to assist with language editing. The authors assume responsibility for all content.

\bibliographystyle{siamplain}
\bibliography{reference}

@article{2017Boundary,
  title = {Boundary First Flattening},
  author = {Sawhney, Rohan and Crane, Keenan},
  journal = {ACM Trans. Graph.},
  volume = {37},
  number = {1},
  year = {2017},
  eid = {5},
}

@article{2008Discrete,
  title = {Discrete Surface {Ricci} Flow},
  author = {Jin, M. and Kim, J. and Luo, F. and Gu, X.},
  journal = {IEEE Trans. Vis. Comput. Graph.},
  volume = {14},
  number = {5},
  pages = {1030--1043},
  year = {2008},
}

@book{GEFW13,
  title = {An Introduction to the Theory of Higher-Dimensional Quasiconformal Mappings},
  author = {Gehring, Frederick W. and Martin, Gaven J. and Palka, Bruce P.},
  series = {Math. Surveys Monogr.},
  volume = {216},
  year = {2017},
  publisher = {American Mathematical Society},
  address = {Providence, RI},
}

@book{IWT01,
  title = {Geometric Function Theory and Non-linear Analysis},
  author = {Iwaniec, Tadeusz and Martin, Gaven},
  year = {2001},
  publisher = {Oxford University Press},
}

@article{choi2018density,
  title = {Density-equalizing maps for simply connected open surfaces},
  author = {Choi, Gary P. T. and Rycroft, Chris H},
  journal = {SIAM J. Imaging Sci.},
  volume = {11},
  number = {2},
  pages = {1134--1178},
  year = {2018},
  publisher = {SIAM},
}

@article{choi2022adaptive,
  title = {Adaptive area-preserving parameterization of open and closed anatomical surfaces},
  author = {Choi, Gary P. T. and Giri, Amita and Kumar, Lalan},
  journal = {Comput. Biol. Med.},
  volume = {148},
  eid = {105715},
  year = {2022},
  publisher = {Elsevier},
}

@article{ABKC18,
  author = {Alex Baden and Keenan Crane and Misha Kazhdan},
  title = {M\"{o}bius Registration},
  journal = {Comput. Graph. Forum},
  volume = {37},
  number = {5},
  pages = {211--220},
  year = {2018},
}

@article{Aldaz2012,
  author = {Aldaz, J. M.},
  title = {Sharp bounds for the difference between the arithmetic and geometric means},
  journal = {Arch. Math. (Basel)},
  volume = {99},
  pages = {393--399},
  year = {2012},
}

@article{CPGP18,
  title = {Efficient feature-based image registration by mapping sparsified surfaces},
  author = {Chun Pang Yung and Gary P. T. Choi and Ke Chen and Lok Ming Lui},
  journal = {J. Vis. Commun. Image Represent.},
  volume = {55},
  pages = {561--571},
  year = {2018},
}

@article{DPGP22,
  author = {Daoping Zhang and Gary P. T. Choi and Jianping Zhang and Lok Ming Lui},
  title = {A Unifying Framework for $n$-Dimensional Quasi-Conformal Mappings},
  journal = {SIAM J. Imaging Sci.},
  volume = {15},
  number = {2},
  pages = {960--988},
  year = {2022},
}

@article{DZKC20,
  title = {3{D} Orientation-Preserving Variational Models for Accurate Image Registration},
  volume = {13},
  number = {3},
  journal = {SIAM J. Imaging Sci.},
  author = {Daoping Zhang and Ke Chen},
  year = {2020},
  pages = {1653--1691},
}

@article{ECKC25,
  title = {Rectangular Surface Parameterization},
  volume = {44},
  eid = {111},
  number = {4},
  journal = {ACM Trans. Graph.},
  author = {Etienne Corman and Keenan Crane},
  year = {2025},
}

@article{JWTL23,
  title = {Ellipsoidal conformal and area-/volume-preserving parameterizations and associated optimal mass transportations},
  volume = {49},
  number = {4},
  eid = {50},
  journal = {Adv. Comput. Math.},
  author = {Jia-Wei Lin and Tiexiang Li and Wen-Wei Lin and Tsung-Ming Huang},
  year = {2023},
}

@article{KCLM14,
  author = {Ka Chun Lam and Lok Ming Lui},
  title = {Landmark- and Intensity-Based Registration with Large Deformations via Quasi-conformal Maps},
  journal = {SIAM J. Imaging Sci.},
  volume = {7},
  number = {4},
  pages = {2364--2392},
  year = {2014},
}

@article{KSNL19,
  title = {Curvature adaptive surface remeshing by sampling normal cycle},
  journal = {Comput.-Aided Des.},
  volume = {111},
  pages = {1--12},
  year = {2019},
  author = {Kehua Su and Na Lei and Wei Chen and Li Cui and Hang Si and Shikui Chen and Xianfeng Gu},
}

@article{KSWC17,
  title = {Volume preserving mesh parameterization based on optimal mass transportation},
  journal = {Comput.-Aided Des.},
  volume = {82},
  pages = {42--56},
  year = {2017},
  author = {Kehua Su and Wei Chen and Na Lei and Junwei Zhang and Kun Qian and Xianfeng Gu},
}

@article{LBPS02,
  author = {Bruno L{\'e}vy and Sylvain Petitjean and Nicolas Ray and J{\'e}rome Maillot},
  title = {Least Squares Conformal Maps for Automatic Texture Atlas Generation},
  journal = {ACM Trans. Graph.},
  year = {2002},
  volume = {21},
  number = {3},
  pages = {362--371},
}

@article{LZXGG08,
  author = {Liu, Ligang and Zhang, Lei and Xu, Yin and Gotsman, Craig and Gortler, Steven J.},
  title = {A Local/Global Approach to Mesh Parameterization},
  journal = {Comput. Graph. Forum},
  volume = {27},
  number = {5},
  pages = {1495--1504},
  year = {2008},
}

@article{LMKC13,
  author = {Lok Ming Lui and Ka Chun Lam and Tsz Wai Wong and Xianfeng Gu},
  title = {Texture Map and Video Compression Using {Beltrami} Representation},
  journal = {SIAM J. Imaging Sci.},
  volume = {6},
  number = {4},
  pages = {1880--1902},
  year = {2013},
}

@article{MDMM02,
  author = {Mathieu Desbrun and Mark Meyer and Pierre Alliez},
  title = {Intrinsic parameterizations of surface meshes},
  journal = {Comput. Graph. Forum},
  year = {2002},
  volume = {21},
  number = {3},
  pages = {209--218},
}

@article{MGBS21,
  year = {2021},
  volume = {40},
  number = {4},
  author = {Mark Gillespie and Boris Springborn and Keenan Crane},
  title = {Discrete conformal equivalence of polyhedral surfaces},
  journal = {ACM Trans. Graph.},
  eid = {103},
}

@article{MHWW17,
  author = {Mei-Heng Yueh and Wen-Wei Lin and Chin-Tien Wu and Shing-Tung Yau},
  title = {An Efficient Energy Minimization for Conformal Parameterizations},
  journal = {J. Sci. Comput.},
  volume = {73},
  number = {1},
  pages = {203--227},
  year = {2017},
}

@article{MHWW19,
  author = {Mei-Heng Yueh and Wen-Wei Lin and Chin-Tien Wu and Shing-Tung Yau},
  title = {A Novel Stretch Energy Minimization Algorithm for Equiareal Parameterizations},
  journal = {J. Sci. Comput.},
  year = {2019},
  volume = {78},
  number = {3},
  pages = {1353--1386},
}

@article{MHTL19,
  author = {Mei-Heng Yueh and Tiexiang Li and Wen-Wei Lin and Shing-Tung Yau},
  title = {A Novel Algorithm for Volume-Preserving Parameterizations of 3-Manifolds},
  journal = {SIAM J. Imaging Sci.},
  volume = {12},
  number = {2},
  pages = {1071--1098},
  year = {2019},
}

@article{MHTM21,
  year = {2021},
  volume = {88},
  number = {3},
  eid = {64},
  author = {Mei-Heng Yueh and Tsung-Ming Huang and Tiexiang Li and Wen-Wei Lin and Shing-Tung Yau},
  title = {Projected Gradient Method Combined with Homotopy Techniques for Volume-Measure-Preserving Optimal Mass Transportation Problems},
  journal = {J. Sci. Comput.},
}

@article{QWWX23,
  title = {Practical construction of globally injective parameterizations with positional constraints},
  volume = {9},
  number = {2},
  journal = {Comput. Vis. Media},
  author = {Qi Wang and Wen-Xiang Zhang and Yuan-Yuan Cheng and Ligang Liu and Xiao-Ming Fu},
  year = {2023},
  pages = {265--277},
}

@article{RLZW14,
  author = {Rongjie Lai and Zaiwen Wen and Wotao Yin and Xianfeng Gu and Lok Ming Lui},
  title = {Folding-Free Global Conformal Mapping for Genus-0 Surfaces by Harmonic Energy Minimization},
  journal = {J. Sci. Comput.},
  year = {2014},
  volume = {58},
  number = {3},
  pages = {705--725},
}

@article{SHSA00,
  author = {Steven Haker and Sigurd Angenent and Allen Tannenbaum and Ron Kikinis and Guillermo Sapiro and Michael Halle},
  title = {Conformal surface parameterization for texture mapping},
  journal = {IEEE Trans. Vis. Comput. Graph.},
  year = {2000},
  volume = {6},
  number = {2},
  pages = {181--189},
}

@article{TMWH23,
  author = {Tsung-Ming Huang and Wei-Hung Liao and Wen-Wei Lin and Mei-Heng Yueh and Shing-Tung Yau},
  title = {Convergence Analysis of Volumetric Stretch Energy Minimization and Its Associated Optimal Mass Transport},
  journal = {SIAM J. Imaging Sci.},
  volume = {16},
  number = {3},
  pages = {1825--1855},
  year = {2023},
}

@article{TMWH24,
  title = {Fundamental theory and {R}-linear convergence of stretch energy minimization for spherical equiareal parameterization},
  author = {Tsung-Ming Huang and Wei-Hung Liao and Wen-Wei Lin},
  pages = {1--25},
  volume = {32},
  number = {1},
  journal = {J. Numer. Math.},
  year = {2024},
}

@inproceedings{WCSS24,
  author = {Wei Chen and Siquan Sun and Yue Wang and Na Lei and Chander Sadasivan and Apostolos Tassiopoulos and Shikui Chen and Hang Si and Xianfeng Gu},
  title = {Robust Surface Remeshing Based on Conformal Welding},
  booktitle = {Proceedings of the 2024 International Meshing Roundtable (IMR)},
  pages = {28--39},
  year = {2024},
  editor = {Ketan Mittal and Jessica Zhang and Scott Mitchell},
  publisher = {SIAM},
}

@article{WWCJ21,
  title = {3{D} brain tumor segmentation using a two-stage optimal mass transport algorithm},
  volume = {11},
  number = {1},
  eid = {14686},
  journal = {Sci. Rep.},
  author = {Wen-Wei Lin and Cheng Juang and Mei-Heng Yueh and Tsung-Ming Huang and Tiexiang Li and Sheng Wang and Shing-Tung Yau},
  year = {2021},
}

@article{WWJW22,
  title = {A novel 2-phase residual {U}-net algorithm combined with optimal mass transportation for 3{D} brain tumor detection and segmentation},
  volume = {12},
  number = {1},
  eid = {6452},
  journal = {Sci. Rep.},
  author = {Wen-Wei Lin and Jia-Wei Lin and Tsung-Ming Huang and Tiexiang Li and Mei-Heng Yueh and Shing-Tung Yau},
  year = {2022},
}

@article{YTKC15,
  title = {Landmark-Matching Transformation with Large Deformation Via $n$-dimensional Quasi-conformal Maps},
  volume = {67},
  number = {3},
  journal = {J. Sci. Comput.},
  author = {Yin Tat Lee and Ka Chun Lam and Lok Ming Lui},
  year = {2016},
  pages = {926--954},
}

@article{ZHTL25S,
  title = {A robust {H}essian-based trust region algorithm for spherical conformal parameterizations},
  volume = {68},
  number = {6},
  journal = {Sci. China Math.},
  author = {Zhong-Heng Tan and Tiexiang Li and Wen-Wei Lin and Shing-Tung Yau},
  year = {2025},
  pages = {1461--1486},
}

@article{ZHTL25,
  title = {$n$-Dimensional Volumetric Stretch Energy Minimization for Volume-/Mass-Preserving Parameterizations with Applications},
  volume = {18},
  number = {2},
  journal = {SIAM J. Imaging Sci.},
  author = {Zhong-Heng Tan and Tiexiang Li and Wen-Wei Lin and Shing-Tung Yau},
  year = {2025},
  pages = {1141--1175},
}

@article{ZZHW25,
  title = {{OMT} and tensor {SVD}-based deep learning model for segmentation and predicting genetic markers of glioma: A multicenter study},
  volume = {122},
  number = {28},
  journal = {Proc. Natl. Acad. Sci. USA},
  author = {Zhengyang Zhu and Han Wang and Tiexiang Li and Tsung-Ming Huang and Huiquan Yang and Zhennan Tao and Zhong-Heng Tan and Jianan Zhou and Sixuan Chen and Meiping Ye and Zhiqiang Zhang and Feng Li and Dongming Liu and Maoxue Wang and Jiaming Lu and Wen Zhang and Xin Li and Qian Chen and Zhuoru Jiang and Futao Chen and Xin Zhang and Wen-Wei Lin and Shing-Tung Yau and Bing Zhang},
  year = {2025},
  eid = {e2500004122},
}

@misc{AIM,
  author = {{AIM@SHAPE}},
  title = {{Digital Shape Workbench - Shape Repository}},
  howpublished = {\url{http://visionair.ge.imati.cnr.it/ontologies/shapes/}},
  year = {2016},
}

@misc{Stanford,
  author = {{Stanford Computer Graphics Laboratory}},
  title = {{The Stanford 3D Scanning Repository}},
  howpublished = {\url{https://graphics.stanford.edu/data/3Dscanrep/}},
  year = {2023},
}

@article{Menze2015BRATS,
  author = {Menze, Bjoern H. and Jakab, Andras and Bauer, Stefan and
             Kalpathy-Cramer, Jayashree and Farahani, Keyvan and
             Kirby, Justin and others},
  title = {The Multimodal Brain Tumor Image Segmentation Benchmark ({BRATS})},
  journal = {IEEE Trans. Med. Imaging},
  volume = {34},
  number = {10},
  pages = {1993--2024},
  year = {2015},
}

@article{Bakas2017TCGA,
  author = {Bakas, Spyridon and Akbari, Hamed and Sotiras, Aristeidis and
             Bilello, Michel and Rozycki, Martin and Kirby, Justin S. and others},
  title = {Advancing {The Cancer Genome Atlas} glioma {MRI} collections
             with expert segmentation labels and radiomic features},
  journal = {Sci. Data},
  volume = {4},
  year = {2017},
  eid = {170117},
}

@misc{Bakas2018BRATS,
  author = {Bakas, Spyridon and Reyes, Mauricio and Jakab, Andras and
             Bauer, Stefan and Rempfler, Markus and Crimi, Alessandro and others},
  title = {Identifying the Best Machine Learning Algorithms for Brain Tumor
             Segmentation, Progression Assessment, and Overall Survival
             Prediction in the {BRATS} Challenge},
  year = {2018},
  howpublished = {arXiv preprint},
  eprint = {1811.02629},
}

@book{Aubin1998,
  author = {Thierry Aubin},
  title = {Some Nonlinear Problems in Riemannian Geometry},
  series = {Springer Monogr. Math.},
  publisher = {Springer},
  address = {Berlin, Heidelberg},
  year = {1998},
}

@article{DyerVegterWintraecken2015,
  author = {Ramsay Dyer and Gert Vegter and Mathijs Wintraecken},
  title = {Riemannian simplices and triangulations},
  journal = {Geom. Dedicata},
  volume = {179},
  pages = {91--138},
  year = {2015},
}

@book{Hatcher2002,
  author    = {Allen Hatcher},
  title     = {Algebraic Topology},
  publisher = {Cambridge University Press},
  year      = {2002}
}

@article{Lipman2014,
  author = {Yaron Lipman},
  title = {Bijective Mappings of Meshes with Boundary and the Degree in Mesh Processing},
  journal = {SIAM J. Imaging Sci.},
  volume = {7},
  number = {2},
  pages = {1263--1283},
  year = {2014},
}

@misc{Crane19,
title = {The $n$-dimensional cotangent formula},
author = {Keenan Crane},
year = {2019},
url = {https://www.cs.cmu.edu/~kmcrane/Projects/Other/nDCotanFormula.pdf},
}

@article{klyh:2021,
  author = {Kuo, Yueh-Cheng and Lin, Wen-Wei and Yueh, Mei-Heng and Yau, Shing-Tung},
  title = {Convergent Conformal Energy Minimization for the Computation of Disk Parameterizations},
  journal = {SIAM J. Imaging Sci.},
  volume = {14},
  number = {4},
  pages = {1790--1815},
  year = {2021},
}

@article{LY25,
  author = {Shu-Yung Liu and Mei-Heng Yueh},
  title = {Energy-Based Distortion-Balancing Parameterization for Open Surfaces},
  journal = {SIAM J. Imaging Sci.},
  volume = {18},
  number = {4},
  pages = {2059--2093},
  year = {2025},
}

@book{NW06,
  author = {Jorge Nocedal and Stephen J. Wright},
  title = {Numerical Optimization},
  edition = {2nd},
  publisher = {Springer},
  address = {New York},
  year = {2006},
}

@book{AMS08,
  author    = {P.-A. Absil and Robert Mahony and Rodolphe Sepulchre},
  title     = {Optimization Algorithms on Matrix Manifolds},
  publisher = {Princeton University Press},
  address   = {Princeton, NJ},
  year      = {2008},
  isbn      = {978-0-691-13298-3}
}

@article{HGA15,
  author = {Wen Huang and Kyle A. Gallivan and P.-A. Absil},
  title = {A {Broyden} Class of Quasi-{Newton} Methods for {Riemannian} Optimization},
  journal = {SIAM J. Optim.},
  volume = {25},
  number = {3},
  pages = {1660--1685},
  year = {2015},
}

@misc{stanford-david,
  author       = {{Stanford Computer Graphics Laboratory}},
  title        = {{Digital Michelangelo Project Archive}},
  howpublished = {\url{https://graphics.stanford.edu/data/mich/}}
}

@misc{vcg-models,
  author       = {{Visual Computing Laboratory, ISTI--CNR}},
  title        = {{VCG Library and Public Mesh-Model Resources}},
  howpublished = {\url{https://vcglib.net/}}
}

@misc{sketchfab,
  author       = {{Sketchfab}},
  title        = {{Public 3D Model Repository}},
  howpublished = {\url{https://sketchfab.com/}}
}

@misc{Keenan,
  author       = {Keenan Crane},
  title        = {{Keenan's 3D Model Repository}},
  howpublished = {\url{https://www.cs.cmu.edu/~kmcrane/Projects/ModelRepository/}}
}
\end{document}